\documentclass[12pt]{amsart}

\usepackage{times}
\usepackage{longtable}
\usepackage{url}
\usepackage{enumerate}
\usepackage{epsfig}

\newcommand{\Ord}{{\mathcal O}}                  
\newcommand{\alg}[1]{\mathcal{#1}}                
\DeclareMathOperator{\supp}{{\mathrm supp}}                   
\DeclareMathOperator{\tr}{{\mathrm tr}}                   
\DeclareMathOperator{\co}{{\mathrm co}}                   
\renewcommand{\Re}{\mathop{\mathrm{Re}}}
\renewcommand{\Im}{\mathop{\mathrm{Im}}}
\newcommand{\ST}{\mathcal{S}}                    

\newcommand{\IC}{\mathbf{C}}                     
\newcommand{\IR}{\mathbf{R}}                     
\newcommand{\IZ}{\mathbf{Z}}                     
\newcommand{\eps}{{\varepsilon}}                 

\theoremstyle{plain} 
\newtheorem{Theorem}{Theorem}[section]
\newtheorem{Lemma}[Theorem]{Lemma}
\newtheorem{Proposition}[Theorem]{Proposition}

\newtheorem{Corollary}[Theorem]{Corollary}

\theoremstyle{definition} 
\newtheorem{Definition}[Theorem]{Definition}

\newtheorem{Example}[Theorem]{Example}
 
\numberwithin{equation}{section}

\begin{document}
\title[Some examples of Brown measures]{Computation of some examples of Brown's spectral measure in free probability}
\author{Philippe Biane}
\address{CNRS, DMA\\ 45 rue d'Ulm\\
75005 Paris\\ France}
\email{Philippe.Biane@ens.fr}
\author{Franz Lehner}
\address{Centre \'Emile Borel\\ Institut Henri Poincar\'e\\
11 rue Pierre et Marie Curie\\
75231 Paris Cedex 05\\
France
}
\email{lehner@bayou.uni-linz.ac.at}
\date{31.12.1999}


\maketitle{}

\begin{abstract}
We use free probability techniques for computing spectra and Brown measures
 of
some non hermitian operators in finite von Neumann algebras. Examples include
$u_n+u_{\infty}$ where $u_n$ and $u_{\infty}$
are the generators of $\IZ_n$ and $\IZ$ respectively, in the free product 
$\IZ_n*\IZ$, or elliptic elements, of the form  $S_{\alpha}+iS_{\beta}$
where $S_\alpha$ and $S_\beta$ are free
semi-circular elements of variance $\alpha$ and $\beta$.
\end{abstract}
\section{Introduction}
\label{sec:introduction}

Recently Haagerup and Larsen 
\cite{HaagerupLarsen:1999:BrownMeasure} have computed the spectrum and 
the Brown measure of $R$-diagonal elements
in a finite von Neumann algebra, in
terms of the  distribution of its radial part. The purpose of this paper is to
apply free probability techniques for computing spectra and Brown measures of
some non-hermitian, and non-$R$-diagonal elements 
in finite von Neumann algebras, which can be written as a free sum of an
$R$-diagonal element and an element with arbitrary $*$-distribution.
Motivations for  this study are twofold, on one hand
some of these elements appear as transition operators
 of random walks on groups or semi-groups, 
 see e.g. \cite{HarpRobVal:1993:spectreI},
\cite{HarpRobVal:1993:spectreII}, \cite{BeVaZuk:1997:Harper}, here we shall for
example treat linear combinations of $u_n$ and $u_{\infty}$,
the generators of $\IZ_n$ and $\IZ$  in  
$\IZ_n*\IZ$ and $u_2+v_2+u_\infty$.
On the other hand random matrix theory has a close connection with 
free probability (see \cite{VoicDykNic:1992:freeRV} ), but for the moment very
little has been done for understanding limit distributions of spectra of 
non-normal matrices in terms of free probability.
For example, the empirical distribution
on the eigenvalues of a random matrix with
independent identically distributed complex entries, suitably rescaled, 
 converges, with probability one, as its size grows to
infinity, to the circular law (the uniform distribution on the unit disk), see 
\cite{Girko:1984:circular}, \cite{Girko:1997:strongcircular},
\cite{Bai:1997:circular}, which
 is the Brown measure for a circular element, in the sense of
Voiculescu. It is known that
the circular element is the limit in $*$-distribution of the above
random matrices,
but it is not possible to deduce from this the convergence
of the empirical distribution on the spectrum (see Lemma 2.1 below).

Another example that we shall consider in this paper is the free sum of an
arbitrary element with a circular element. Hopefully, the corresponding Brown
measures should represent limit of eigenvalue distributions of random matrices
of the form $A+W$ where $A$ is a matrix with some limit $*$-distribution, and
$W$ is a matrix with independent entries.
In addition to the circular element discussed above,
this is known to be true for the so-called \emph{elliptic element},
which can be written as $S_\alpha+i S_\beta$ and whose Brown measure was first
computed in \cite{Larsen:1999:Thesis} by ad-hoc methods.
It turns out to be treatable by our method as well.
The empiciral eigenvalue distribution of its matrix model
with Gaussian random matrices is computed in \cite{HiaiPetz:1998:LogarithmicEnergy}
and shown to converge to the uniform measure on its spectrum, an ellipse.

However in this paper we shall  stick to
the purely free probabilistic aspects of the subject, and not touch upon the
random matrix problem.
We hope to deal with this somewhere else. 

This paper is organized as follows. In section 2 we recall preliminary facts
about Brown measures and free probability theory. In section 3 we give a
general approach towards the computation of the Brown measure for  the sum of
an $R$-diagonal element with an arbitrary element. We specialize in sections 3
and 4 to the cases where the $R$-diagonal element is a Haar unitary or a
circular element, respectively.  We close with some final remarks in section 5.
The pictures of random matrix spectra appearing in various sections of this papers
were computed with \url{GNU} \url{octave} and plotted with \url{gnuplot};
the plots of densities of various Brown measures, which accompany or
replace the rather unwieldy density formulae, were computed by \url{Mathematica}.

\emph{Acknowledgements.}
This work was started during a special semester at the Erwin Schr\"odinger
Institute in Vienna in spring 1999, organized by the Institut f\"ur Funktionalanalysis
of the university of Linz and its head J.B.~Cooper.
The second author was supported by the EU-network ``Non-commutative geometry'' ERB FMRX CT960073.

\section{Preliminaries}
\label{sec:preliminaries}

\subsection{The Fuglede--Kadison determinant and Brown's spectral measure}
Let $\alg M$ be a finite von Neumann algebra with faithful tracial state $\tau$
and denote, for invertible $a\in\alg M$,  $\Delta(a)=e^{\tau(\log|a|)}$ 
its Fuglede-Kadison determinant (cf.\ \cite{FugledeKadison:1952:Determinant}).
Denoting by $\mu_x$ the spectral measure for the self-adjoint element $x\in \alg M$,
i.e.\ the unique probability measure on the real line satisfying
$\tau(x^n) = \int t^n d\mu(t)$, we have the following formula for the 
logarithm of the determinant, which serves as a definition of the determinant
in the case where $a$ is not invertible
$$
\log\Delta(a) = \int_\IR \log t \, d\mu_{|a|}(t)
.
$$
The function $\Delta(\lambda -a)$ is a subharmonic function of the
complex variable $\lambda$, and
there is a unique probability measure $\mu_a$ 
on $\IC$, with support on the spectrum of
$a$, called the 
\emph{Brown measure} of $a$, such that 
$$
\log\Delta(\lambda - a)=\int\log\vert\lambda-z\vert\,\mu_a(dz)
;
$$
it is given by
$$
\mu_a = \frac{1}{2\pi}
        \nabla^2
        \log\Delta(\lambda  - a)
$$
where $\nabla^2$ is the Laplace operator in the complex plane,
in the sense of distributions (see \cite{Brown:1986:BrownMeasure}).
If $a$ is normal, then $\mu_a$ is just the spectral measure of $a$. When
$\alg M$ is $M_n(\IC)$, with the canonical normalized trace,
 then $\mu_a$ is the empirical distribution on the
spectrum of $a$ (counting multiplicities).
Although the Brown measure of $a$ 
can be computed from its $*$-distribution,
i.e. the collection of all its $*$-moments 
$\tau(a^{\eps_1}a^{\eps_2}\cdots a^{\eps_n})$, where
$a^{\eps_j}$ is either $a$ or $a^*$, it does not depend continuously on these
$*$-moments. Indeed let for example
$a_n$ be the $n\times n$ nilpotent matrix with ones on the first 
upper diagonal and zeros everywhere else, then as $n$ goes to infinity
the $*$-moments of $a_n$
converge towards those of a Haar unitary (a unitary element $u$ with
$\tau(u^n)=0$ for $n\ne0$), whose Brown measure is the Haar measure on the
unit circle,  whereas the Brown measure of $a_n$ is $\delta_0$
for all $n$.
\begin{Lemma}
  \label{lmm:balayage}
  Let $(a_n;n\geq 0)$ be a uniformly bounded
  sequence  whose $*$-distributions
  converge towards that of $a$, and suppose
  the Brown measure of $a_n$ converges weakly
  towards some measure $\mu$, then one has
  \begin{enumerate}[(i)]
   \item $\int\log\vert\lambda-z\vert \,\mu(dz)
    \leq\Delta(\lambda-a)=\int\log\vert\lambda-z\vert \,\mu_a(dz)$
    for all $\lambda\in\IC$ 
   \item $\int\log\vert\lambda-z\vert \,\mu(dz)=\Delta(\lambda-a)
   =\int\log\vert\lambda-z\vert \,\mu_a(dz)$ 
   for all $\lambda$ large enough.
  \end{enumerate}
\end{Lemma}
\begin{proof}
  The distribution of $\vert\lambda-a_n\vert$
  has a support which remains in a fixed
  compact set, 
  and it converges weakly towards that of $\vert\lambda-a\vert$.
  Part $(i)$ follows from this and 
  the fact that the function $\log$ is a limit of a
  decreasing sequence of continuous functions. If $\lambda$ is large enough, then 
  the union of the supports of the distributions of the $\vert \lambda-a_n\vert$
  is away from $0$, hence the function $\log$ is continuous there and $(ii)$
  follows from weak convergence.
\end{proof}

The outcome of $(i)$ of
the Lemma is that the measure $\mu_a$ is a balay\'ee of measure
 $\mu$, while we get from $(ii)$ the following
\begin{Corollary}
  Let $U_a$ be the  
  unbounded connected component of the complement of the support of $\mu_a$, 
  then the support of $\mu$ is included in $\IC\setminus U_a$.
\end{Corollary}
\begin{proof}
  The function $\int\log\vert\lambda-z\vert\,\mu_a(dz)$ is harmonic in
  $U_a$, while $\int\log\vert\lambda-z\vert\,\mu(dz)$ is subharmonic there,
  consequently 
  $
  \int\log\vert\lambda-z\vert\,\mu_a(dz)
  -\int\log\vert\lambda-z\vert\,\mu(dz)
  $ 
  is a nonnegative superharmonic
  function on $U_a$. Since this function attains the value $0$ by $(ii)$, it is
  identically $0$ by the minimum principle, therefore 
  $\int\log\vert\lambda-z\vert\,\mu(dz)$ is harmonic on $U_a$, and thus the support
  of $\mu$ is included in $\IC\setminus U_a$.
\end{proof}
Conversely, given  two measures $\mu$ and $\mu_a$
on $\IC$ satisfying $(i)$ and $(ii)$, we do not know
whether there always exists a corresponding sequence $(a_n)_{n\geq0}$,
fulfilling the hypotheses of Lemma~\ref{lmm:balayage}

\subsection{$R$- and $S$-transforms}
We shall refer to \cite{VoicDykNic:1992:freeRV},
and  \cite{Voiculescu:1999:SaintFlour} or \cite{HiaiPetz:1999:book}
for basic concepts of free probability theory.
Let $(\alg M, \tau)$ be as in section 2.1, and let $a\in \alg M$. 
The power series
$$
G_a(\zeta)
 = \frac{1}{\zeta}
   \sum_{n=0}^\infty \frac{\tau(a^n)}{\zeta^n}
$$
can be inverted (for composition of formal power series), in the form
$$
K_a(z) = \frac{1}{z} + \sum_{n=0}^\infty c_{n+1} z^n
         = \frac{1}{z}
           (1+R_a(z))
.
$$
The power series $R_a$ is called the \emph{$R$-transform} of $a$
(note that this slightly differs from the original definition of Voiculescu)
and its coefficients are called the
\emph{free cumulants} of $a$.
Let
$$
\psi_a(z) = \sum_{n=1}^\infty \tau(a^n) \,z^n
 = \frac{1}{z}\,
   G_a\left(
        \frac{1}{z}
      \right)
   -1
$$
be the generating moment series for $a$, and 
assume that the first moment is nonzero, so that $\psi_a'(0)\ne 0$.
Then $\psi_a$ has an inverse $\chi_a$, and the $S$-transform
of $a$ is defined as
$$
S_a(z)=\frac{1+z}{z}\, \chi_a(z)
$$ 
Observe that the power series $zS_a(z)$ and
 $R_a(z)$ are then inverse of each other (when the mean is nonzero).
The relevance of these series to free probability is that, if $a,b\in \alg M$
are free, then
$$R_{a+b}=R_a+R_b\qquad\text{and}\qquad S_{ab}=S_aS_b$$
see e.g. \cite{VoicDykNic:1992:freeRV}.
\subsection{Calculus of $R$-diagonal elements}
\label{subsec:RDiagonalCalculus}
We use the same notations as in the previous section.
\begin{Definition}
  A non-commutative random variable $x$ is called \emph{$R$-diagonal}, if
  $x$ has polar decomposition $x=uh$, where $u$ is a Haar unitary
  free from the radial part $h=|x|$.
\end{Definition}
 Recall that a unitary $u\in\alg M$ is
  called a Haar unitary if $\tau(u^n)=0$ for all integers $n\ne0$.
One can check that 
the product of an arbitrary element $y$ with a free Haar unitary is an
$R$-diagonal element.
According to 
\cite{HaagerupLarsen:1999:BrownMeasure},
any $R$-diagonal element with polar decomposition
$x=uh$ has the same distribution as a product
  $a\tilde{h}$,
where $\tilde{h}$ has a symmetric distribution, and its absolute value is
distributed as $h$, whereas  $a$ is a self-adjoint unitary,
free from $\tilde{h}$, and of zero trace.
Indeed, one can assume 
$\tilde{h}=a'h$, where $a'$ is a symmetry commuting with $h$ and 
$aa'$ is a Haar unitary free from $h$.
Let $a,b$ be two free $R$-diagonal elements, then one has
equality in $*$-distribution of the pairs  $(a,b)$ and $(ua,ub)$
 where $u$ is a Haar
unitary free with $\{a,b\}$, therefore $a+b$ has the same $*$-distribution
 as $u(a+b)$ which is $R$-diagonal, and thus the sum of two free
$R$-diagonal elements is again $R$-diagonal. 
Let $f_x(z^2)=R_{\tilde h}(z)$
 be the cumulant series of $\tilde{h}$, which determines the
$*$-distribution of $x$, then
the power series $z(1+z)S_{x^*x}(z)$ and $f_x(z)$ are inverse of each other.
Furthermore
if $a,b$ are two free $R$-diagonal elements, then one has
\begin{equation}
  \label{eq:sumofrdiag}
  f_{a+b}=f_a+f_b.
\end{equation}
See \cite{NicaSpeicher:1997:Rdiagonal}, \cite{NicaSpeicher:1998:Commutator}
and \cite{HaagerupLarsen:1999:BrownMeasure}.

\subsection{Brown measure of $R$-diagonal elements}

In \cite{HaagerupLarsen:1999:BrownMeasure} the Brown measure of an
$R$-diagonal element is determined as follows.

\begin{Theorem}[{\cite[Thm.~4.4, Prop.~4.6]
{HaagerupLarsen:1999:BrownMeasure}}]
  \label{thm:HaagerupLarsen:5.1}
  Let $u$, $h$ be $*$-free random variables in $(\alg M,\tau)$,
  with $u$ a Haar unitary and $h$ positive s.t.\ the distribution $\mu_h$
  of $h$ is not a Dirac measure.
  Then the Brown measure $\mu_{uh}$ of $uh$ has the following properties.
  \begin{enumerate}[$(i)$]
   \item $\mu_{uh}$ is rotation invariant and its support is the annulus
    with inner radius $\|h^{-1}\|_2^{-1}$ and outer radius $\|h\|_2$.
   \item The $\ST$-transform $\ST_{\mu_{h^2}}$ of $h^2$
    has an analytic continuation to a neighbourhood of\/ \mbox{$]\mu_h(\{0\})-1,0]$}
    and its derivative $\ST_{\mu_{h^2}}'$ is strictly negative on this interval and
    its range is
    $\ST_{\mu_{h^2}}({}]\mu_h(\{0\})-1,0]) = [\|h\|_2^{-2}, \|h^{-1}\|_2^2[$.
   \item $\mu_{uh}(\{0\})=\mu_h(\{0\})$ and for $t\in{}]\mu_h(\{0\}),1]$
    $$
    \mu_{uh}\left(
              B\left(
                 0,\frac{1}{\sqrt{\ST_{\mu_{h^2}}(t-1)}}
               \right)
            \right)
    = t
    $$
   \item $\mu_{uh}$ is the only rotation symmetric probability measure satisfying
    (iii).
   \item If $h$ is invertible then $\sigma(uh)=\supp\mu_{uh}$,
    i.e., the annulus discussed above.
   \item \label{item:HaagerupLarsen:5.1:notinvertible}
    If $h$ is not invertible then $\sigma(uh)=B(0,\|h\|_2)$.
  \end{enumerate}
\end{Theorem}

The proof involves a formula for the spectral radius of products of free elements.
\begin{Proposition}[{\cite[Prop.~4.1]{HaagerupLarsen:1999:BrownMeasure}}]
  \label{prop:HaagerupLarsen:4.10}
  Let $a,b$ be $*$-free centered elements in $\alg M$.
  Then the spectral radius of $ab$ is
  $$
  \rho(ab) = \|a\|_2 \, \|b\|_2
  $$
\end{Proposition}

In particular, an $R$-diagonal element $a=uh$ can be written as
$u_1u_2h$, with free Haar unitaries $u_1$, $u_2$ and therefore its spectral radius is
$\rho(a)=\|u_1\|_2\,\|u_2h\|_2 = \|a\|_2$.

\section{Adding an $R$-diagonal element}
In this section we give a general approach to computing the
Brown measure of the
  sum
of a random variable with an arbitrary distribution and  a free
 $R$-diagonal element.
So we let $a$ be an arbitrary element, 
$h$ be self-adjoint and $u$ a Haar unitary, with $\{a,u,h\}$
forming a free family.

\subsection{The spectrum of $a+uh$}

The spectrum of $a+uh$ is determined as follows.
For $\lambda\not\in\sigma(a)$, 
$\lambda - a - uh$ is invertible if and only if $1-uh(\lambda-a)^{-1}$
is invertible.
If $h$ is not invertible, then by the result of Haagerup and Larsen on 
$R$-diagonal elements, the latter is the case
if and only if 
\begin{equation}
  \label{eq:spectrumupperbound}
  \|h(\lambda-a)^{-1}\|_2 = \|h\|_2 \, \|(\lambda-a)^{-1}\|_2 < 1
;
\end{equation}
if $h$ is invertible, we get the additional possibility that
$1<\|h^{-1}\|_2\, \|\lambda-a\|_2$.
In this case we can look at $(uh)^{-1}(\lambda-a)-1$.

The case where $\lambda\in\sigma(a)$ must be considered individually.
Complications arise for such $\lambda$, for which $\lambda\in\sigma(a)$,
but $\|(\lambda-a)^{-1}\|_2<\infty$.
Otherwise condition \eqref{eq:spectrumupperbound}
will be satisfied when approaching $\lambda$ from the outside of $\sigma(a)$,
so that $\lambda$ lies in the closure of the spectrum of $a+uh$,
hence in the spectrum.

\subsection{The Brown measure of $a+uh$}
\label{sec:a+uhBrown}
We can assume that  $u=u_1^*u_2$ with $u_1$ and $u_2$ Haar unitaries, where 
$\{u_1,u_2,a,h\}$ is a free family, to get
\begin{align*}
  \log\Delta(\lambda-a-uh)
  &= \tau(\log|u_1^*(u_1(\lambda-a)-u_2h)|)\\
  &= \tau(\log|u_1(\lambda-a)-u_2h|)\\
  &= \int \log|z|\,d\mu_{u_1(\lambda-a)-u_2h}(z)
\end{align*}

and this is the Fuglede--Kadison determinant of
$x_\lambda=u_1(\lambda-a)-u_2h$, which is an $R$-diagonal element whose
$*$-distribution can be computed according to \eqref{eq:sumofrdiag}, i.e. 
$f_x = f_{u_1|\lambda-a|}+f_{u_2h}$. This in turn will yield the $S$-transform
of  $x_\lambda^*x_\lambda$, and then by
 Theorem~\ref{thm:HaagerupLarsen:5.1}, we can compute
 $ \log\Delta(\lambda-a-uh)$.

From the discussion in section~\ref{subsec:RDiagonalCalculus}
we have the relation
\begin{equation}
  \label{eq:fxl-1=Rxl2}
f_{x_\lambda}^{\langle-1\rangle}(\zeta)
 = \frac{1}{\zeta}
   \left(
     1
     +
     R_{x_\lambda^*x_\lambda}\left(
                               \frac{1}{\zeta}
                             \right)
   \right)
\end{equation}
In order to be more specific, let us assume that
 $a$ is self-adjoint, then 
 the computation of the distribution of $(\lambda-a)^*(\lambda-a)$
is conveniently accomplished by using the Cauchy transform of $a$, 
namely factoring $\zeta-|\lambda-x|^2 = (x-x_+)(x-x_-)$ with
\begin{equation}
  \label{eq:Gl-a2:xplusminus}
  x_{\pm}
  = \frac{1}{2}
    \left(
      \lambda + \bar{\lambda}
      \pm
      \sqrt{(\lambda-\bar{\lambda})^2+4\zeta}
    \right)
  = \Re\lambda \pm i\sqrt{(\Im\lambda)^2-\zeta}
\end{equation}
and expanding into partial fractions
$$
\frac{1}{\zeta-|\lambda-x|^2}
= \frac{1}{x_+ - x_-}
  \left(
    \frac{1}{x_+ - x}
    -
    \frac{1}{x_- - x}
  \right)
$$
we get
\begin{equation}
  \label{eq:G|l-a|^2}
  G_{|\lambda-a|^2}(\zeta)
  = \int \frac{d\mu_a(x)}%
              {\zeta-|\lambda-x|^2}
  = \frac{G_a(x_+)- G_a(x_-)}%
         {x_+ - x_-}
 .
\end{equation}
Using the same technique one can compute the $2$-norm of the inverse of $\lambda-a$.
Assuming again that $a$ is self-adjoint we have that
\begin{equation}
  \label{eq:||(l-a)^(-1)||_2}
  \begin{aligned}
    \|(\lambda-a)^{-1}\|_2^2
    &= \int\frac{d\mu_a(x)}{|\lambda-x|^2}\\
    &= \int\frac{d\mu_a(x)}{(\lambda-x)(\bar{\lambda}-x)}\\
    &= \frac{1}{\lambda-\bar{\lambda}}
       \int
        \left(
          \frac{1}{\bar{\lambda}-x}
          -
          \frac{1}{\lambda-x}
        \right)
        d\mu_a(x) \\
    &= -\frac{G_a(\lambda)-G_a(\bar{\lambda})}{\lambda-\bar{\lambda}}
  \end{aligned}
\end{equation}

Let us consider the simplest non-trivial random variable,
namely
$
a=u_2=
\left[
  \begin{smallmatrix}
    0&1\\
    1&0
  \end{smallmatrix}
\right]
$ having $2$-point spectrum, so that $|\lambda-u_2|^2$
has a Bernoulli distribution.
The $R$-transform of $|\lambda-u_2|^2=1+|\lambda|^2+ (\lambda+\bar{\lambda})u_2$
is easily computed to be
$$
R_{x_\lambda^*x_\lambda}(z)
= (1+|\lambda|^2)\,z
  +
  \frac{1}{2}
  \left(
    \sqrt{1+4(\lambda+\bar\lambda)^2z^2}
    -
    1
  \right)
$$
and inverting it according to \eqref{eq:fxl-1=Rxl2} leads to an equation of fourth degree,
which apparently is unsuitable for further computations.
So even this simple case seems to be untractable by this method.
In fact, so far we have no concrete example where the general method above can
be carried to the end.
We shall develop other methods, in the next two sections, in order to treat the
cases where the $R$-diagonal element is a Haar  unitary, and a circular element.

\section{Haar unitary case}
\label{sec:Haar}

Now $a$ is an element with an arbitrary distribution, free with a Haar
unitary $u$.

\subsection{The spectrum}
The spectrum of $a+u$ is determined as follows: one has
$\lambda\in\sigma(a+u)$ if and only if 
$1\in\sigma(u^*(\lambda-a))$
and since the latter is $R$-diagonal, we infer from Theorem~\ref{thm:HaagerupLarsen:5.1}
that a necessary and sufficient condition is
\begin{equation}
  \label{eq:sigma(aplusu)}
  \|(\lambda-a)^{-1}\|_2^{-1} \le 1 \le \|(\lambda-a)\|_2  
\end{equation}
if $\lambda\not\in\sigma(a)$; otherwise the condition is simply $1\le\|\lambda-a\|_2$.

\subsection{First approach to the Fuglede-Kadison determinant}
We get the following formula for the Fuglede--Kadison determinant
\begin{equation}
  \label{eq:Haar:logD:integral}
  \begin{aligned}
    \log\Delta(\lambda-a-u)
    &= \tau(\log|\lambda-a-u|)\\
    &= \tau(\log|u^*(\lambda-a)-1|)\\
    &= \int \log|z-1|\, d\mu_{u^*(\lambda-a)}(z)
    .
  \end{aligned}
\end{equation}
Observe that $u^*(\lambda-a)$ is an $R$-diagonal element,
and we can evaluate the integral as follows.
The Brown measure of  an $R$-diagonal element $uh$
  is rotationally symmetric
with radial distribution $\nu(dr)$ and one has
\begin{align*}
  \int \log|z-1|\, d\mu_{uh}(z)
  &= \int_{\|h^{-1}\|_2^{-1}}^{\|h\|_2}
     \int_{0}^{2\pi}
     \log|re^{i\theta}-1|\,d\theta
     \,\nu(dr)
\end{align*}
where the inner integral reduces to
$$
\frac{1}{2\pi}
\int_0^{2\pi}
 \log|re^{i\theta}-1|\,d\theta
= \begin{cases}
    0 & r<1, \\
    \log r
    & r\ge1
    .
  \end{cases}
$$
Introduce the radial distribution function
$$
F_{uh}(r) = \mu_{uh}(B(0,r)) = 
2\pi \int_{\|h^{-1}\|_2^{-1}}^r \nu(d\rho)
$$
which according to 
Theorem~\ref{thm:HaagerupLarsen:5.1}
is related to the moment generating function $\psi_{h^2}$ by
$$
\psi_{h^2}
\left(
  {F_{uh}(r)-1\over F_{uh}(r) \, r^2}
\right)
=F_{uh}(r)-1
$$
(for $\|h^{-1}\|_2^{-1}\le r \le \|h\|_2$),
and by partial integration (note that $F(\|h\|_2)=1$)
\begin{align*}
  \tau(\log|uh-1|)
  &= \int_{\max(1,\|h^{-1}\|_2^{-1})}^{\|h\|_2}
      2\pi \log(r) \, \nu(dr)\\
  &= \log r\, \Bigl. F_{uh}(r)\Bigr|_{\max(1,\|h^{-1}\|_2^{-1})}^{\|h\|_2}
     - \int_{\max(1,\|h^{-1}\|_2^{-1})}^{\|h\|_2} 
     \frac{F_{uh}(\rho)}{\rho}\,d\rho \\
  &= \int_{\max(1,\|h^{-1}\|_2^{-1})}^{\|h\|_2} 
  \frac{1-F_{uh}(\rho)}{\rho}\,d\rho
\end{align*}

\begin{Example}[$2\times2$ matrix]
  Let $a$ have the $*$-distribution of a $2\times 2$ matrix,
   and consider $a+u$, $u$ a Haar unitary.
  Let $\mu_{\pm}=\mu_{\pm}(\lambda)$ be the eigenvalues of
  $|\lambda-a|^2$
  and let
  $$
  G_{|\lambda-a|^2}(\zeta)
  = \frac{1}{2}
    \left(
      \frac{1}{\zeta-\mu_+}
      +
      \frac{1}{\zeta-\mu_-}
    \right)
  $$
  be its Cauchy transform.
  Then 
  $$
  \psi(z) = \frac{1}{z}\,
            G\left(
               \frac{1}{z}
             \right) - 1
          = \frac{1}{2}
            \left(
              \frac{1}{1 - \mu_+ z}
              +
              \frac{1}{1 - \mu_- z}
            \right)
            -1
  $$
  and we get $F(r)$ by solving the equation $\psi(\frac{t-1}{t\, r^2})=t-1$
  for $t$:
  $$
  \left(
    \frac{1}{1 - \mu_+ \frac{t-1}{t\, r^2}}
    +
    \frac{1}{1 - \mu_- \frac{t-1}{t\, r^2}}
  \right)
  = 2t
  $$
  The obvious solution $t=1$ is not interesting for us,
  and dividing it out leads to the other solution
  $$
  F(r) = \frac{2\mu_+\mu_- - r^2(\mu_+\mu_-)}%
              {2(r^2-\mu_+)(r^2-\mu_-)}
       = \frac{\det|\lambda-a|^2 - r^2\tau(|\lambda-a|^2)}%
              {\det(r^2-|\lambda-a|^2)}
  $$
  The logarithm of the Fuglede--Kadison determinant is, for $\lambda\in\sigma(a+u)$,
  \begin{align*}
    \tau(\log|\lambda-a-u|)
    &= \int_1^{\|\lambda-a\|_2}
        \frac{1-F(r)}%
             {r}
        \,
        dr \\
    &= \int_1^{\|\lambda-a\|_2}
        \frac{1}{2}
        \left(
          \frac{r}{r^2-\mu_+}
          +
          \frac{r}{r^2-\mu_-}
        \right)
        dr \\
    &= \left.
         \frac{1}{4}(\log|r^2-\mu_+| + \log|r^2-\mu_-|)
       \right|_1^{\|\lambda-a\|_2} \\
    &= \frac{1}{4}
       \left(
         \log\left|
               \|\lambda-a\|_2^4
               -
               \det|\lambda-a|^2
             \right|
         -
         \log\left|
               1 - 2\|\lambda-a\|_2^2 + \det|\lambda-a|^2
             \right|
       \right)\\
    &= \frac{1}{2}
       \log\left|
             \frac{\mu_+-\mu_-}{2}
           \right|
       -
       \frac{1}{4}
       \left(
         \log|1-\mu_+|
         +
         \log|1-\mu_-|
       \right)
  \end{align*}
  It is now convenient to use the representation of the Laplacian in terms of
  $\partial_\lambda
   = \frac{1}{2}
   \left(
     \frac{\partial}{\partial \Re\lambda}
     -
     i
     \frac{\partial}{\partial \Im\lambda}
   \right)
  $ and its adjoint, namely
  $$
  \nabla^2 = 
     \frac{\partial^2}{\partial (\Re\lambda)^2}
     +
     \frac{\partial^2}{\partial (\Im\lambda)^2}
     = 4\partial_{\bar{\lambda}}\partial_\lambda
  .
  $$
  Then we have the formulae
  \begin{align*}
    \partial_\lambda \|\lambda-a\|_2^2
    &= \partial_\lambda\tau((\lambda-a)^*(\lambda-a)) \\
    &= \tau(\bar{\lambda}-a^*) \\
    \partial_\lambda\det(\lambda-a)
    &= \partial_\lambda( (\lambda-\lambda_1(a))(\lambda-\lambda_2(a)) ) \\
    &= 2\lambda - \lambda_1(a)- \lambda_2(a) \\
    &= 2\tau(\lambda-a)
  \end{align*}
  and the density of the Brown measure of $a+u$ is
  \begin{equation}
    \label{eq:a+u_Density}
    \begin{aligned}
      p_{a+u}(\lambda) 
      &= \frac{2}{\pi}\,
         \partial_{\bar{\lambda}}\partial_\lambda
         \left(
           \log\left|
                 \|\lambda-a\|_2^4
                 -
                 \det|\lambda-a|^2
               \right|
           -
           \log\left|
                 1 - 2\|\lambda-a\|_2^2 + \det|\lambda-a|^2
               \right|
         \right)\\
      &= \frac{2}{\pi}
         \partial_{\bar{\lambda}}
         \left(
           \frac{2\|\lambda-a\|_2^2 \,\tau(\bar{\lambda}-a^*)
                 -2\tau(\lambda-a)\,\det(\bar{\lambda}-a^*) }%
                {\|\lambda-a\|_2^4-\det|\lambda-a|^2} 
  \right.\\& \phantom{=((\frac{1}{2\pi}\partial_{\bar{\lambda}}}\left.\hfill
           -
           \frac{-2\tau(\bar{\lambda}-a^*)
                 +2\tau(\lambda-a)\,\det(\bar{\lambda}-a^*)  }%
                {1-2\|\lambda-a\|_2^2+\det|\lambda-a|^2 }
         \right) \\
      &= \frac{4}{\pi}
         \left(
           \frac{\|\lambda-a\|_2^2 - |\tau(\lambda-a)|^2}%
                { \|\lambda-a\|_2^4-\det|\lambda-a|^2 }
           -2
           \frac{\left|
                   \|\lambda-a\|_2^2\,\tau(\bar{\lambda}-a^*)
                   -\det(\bar{\lambda}-a^*)\,\tau(\lambda-a)
                 \right|^2}%
               {\left(
                   \|\lambda-a\|_2^4-\det|\lambda-a|^2
                \right)^2}
  \right.\\&\phantom{=((((}\left.
           -
           \frac{2|\tau(\lambda-a)|^2-1}%
                {1-2\|\lambda-a\|_2^2+\det|\lambda-a|^2}
           +2
           \frac{\left|
                   \tau(\bar{\lambda}-a^*)
                   -\tau(\lambda-a)\,\det(\bar{\lambda}-a^*)
                 \right|^2}%
                {\left(
                   1-2\|\lambda-a\|_2^2+\det|\lambda-a|^2
                 \right)^2}
         \right)
    \end{aligned}
  \end{equation}
  and in terms of eigenvalues
  \begin{equation}
    \label{eq:a+u_Density_Eigenvalues}
    \begin{aligned}
      p_{a+u}(\lambda)    
      &= \frac{2}{\pi}
         \,
         \partial_{\bar{\lambda}}\partial_\lambda
         \left(
           \frac{1}{2}
           \log\left|
                 \frac{\mu_+-\mu_-}{2}
               \right|
           -
           \frac{1}{4}
           \left(
             \log|1-\mu_+|
             +
             \log|1-\mu_-|
           \right)
         \right)\\
      &= \frac{1}{\pi}
         \,
         \partial_{\bar{\lambda}}
         \left(
           \frac{1}{\mu_+-\mu_-}
           \,
           \partial_\lambda(\mu_+-\mu_-)
           +
           \frac{1}{2}
           \left(
             \frac{1}{1-\mu_+}
             \,
             \partial_\lambda\mu_+
             +
             \frac{1}{1-\mu_-}
             \,
             \partial_\lambda\mu_-
           \right)
         \right) \\
      &= \frac{1}{\pi}
         \left(
           \frac{\partial_{\bar{\lambda}} 
                 \partial_\lambda
                 (\mu_+-\mu_-)}%
                {\mu_+-\mu_-}
           -
           \left|
             \frac{\partial_\lambda
                   (\mu_+-\mu_-)}%
                  {\mu_+-\mu_-}
           \right|^2
  \right.\\&\hskip5em\left.
           +
           \frac{1}{2}
           \left(
             \frac{\partial_{\bar{\lambda}}\partial_\lambda\mu_+}%
                  {1-\mu_+}
             +
             \frac{\partial_{\bar{\lambda}}\partial_\lambda\mu_-}%
                  {1-\mu_-}
             +
             \left|
               \frac{\partial_\lambda\mu_+}%
                    {1-\mu_+}
             \right|^2
             +
             \left|
               \frac{\partial_\lambda\mu_-}%
                    {1-\mu_-}
             \right|^2
           \right)
         \right)
     \end{aligned}
  \end{equation}
  In particular, if
  $
  a=\left[
      \begin{smallmatrix}
        \alpha&0\\
        0&\beta
    \end{smallmatrix}
  \right]
  $
  (Bernoulli distribution) one gets
  $\mu_\pm = \{|\lambda-\alpha|^2|, |\lambda-\beta|^2\}$
  and consequently the density is
  $$
  p_{a+u}(\lambda)  
  = -\frac{|\beta-\alpha|^2}%
          {\pi\left(
             |\lambda-\alpha|^2-|\lambda-\beta|^2   
           \right)^2}
    +
    \frac{1}{2\pi}
    \left(
      \frac{1}{(1-|\lambda-\alpha|^2)^2}
      +
      \frac{1}{(1-|\lambda-\beta|^2)^2}
    \right)
  $$
  on the spectrum, which is determined by the inequalities
  \begin{equation}
    \label{eq:a+u_SpectrumEquation}
    \frac{1}{\mu_+}
    +
    \frac{1}{\mu_-}
    \ge 2
    \qquad
    \mu_+ + \mu_-\ge 2
  \end{equation}

  Specifying further $\alpha=1$, $\beta=-1$, so that $a$ is a symmetry,
  the spectrum is the region bounded by the lemniscate-like curve 
  in the complex plane with the  equation
  $$
  \vert\lambda\vert^2+1=\vert\lambda^2-1\vert^2
  $$
  and we get the picture shown in figure~\ref{fig:u2+uinfty:density}.

  \begin{figure}[pt]
    \begin{center} %
     \includegraphics[width=\textwidth, height=.44\textheight, keepaspectratio]{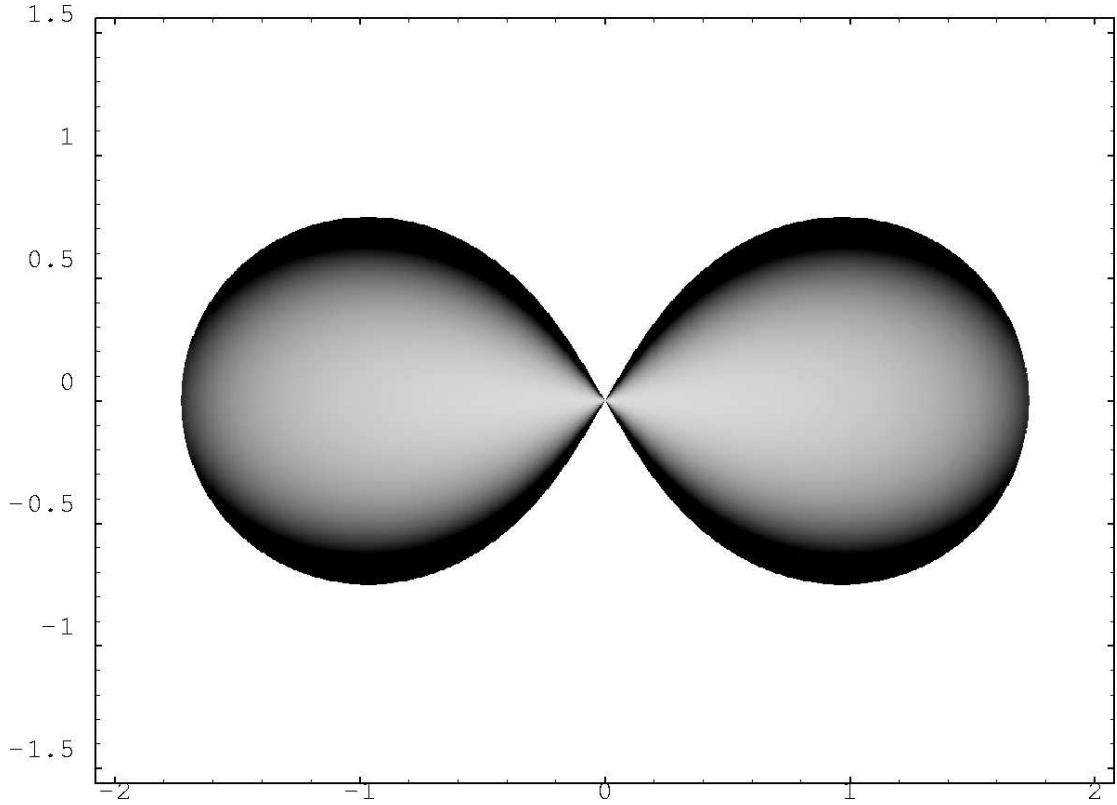}
      \caption{Density of $\mu_{u_2+u_\infty}$}
      \label{fig:u2+uinfty:density}
    \end{center}
  \end{figure}
  \begin{figure}[pb]
    \begin{center} %
     \includegraphics[width=\textwidth, height=.44\textheight, keepaspectratio]{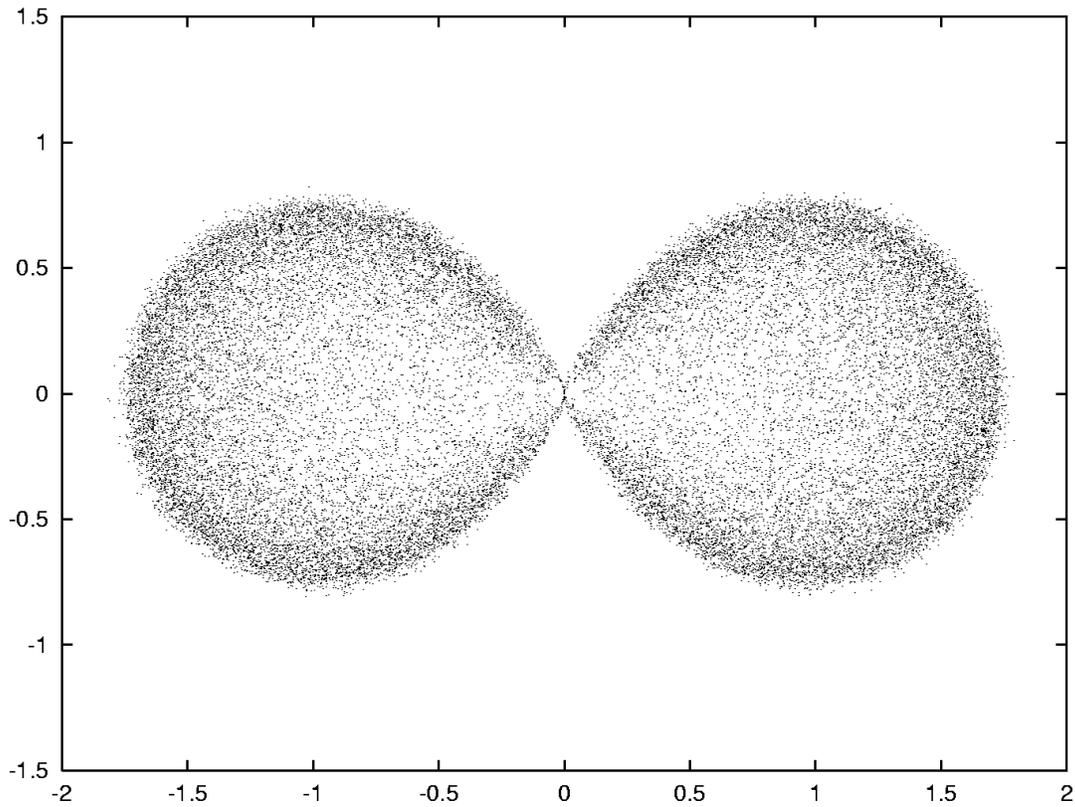}
      \caption{$200$ samples of eigenvalues of $150\times150$ random matrices $U_2+U_\infty$}
      \label{fig:U2+Uinfty_150}
    \end{center}
  \end{figure}
  This should be compared with the sample fig.~\ref{fig:U2+Uinfty_150}
  of eigenvalues of random
  $2N\times 2N$ matrices of the form
  $X=U_2+U_{\infty}$, where $U_{\infty}$ is chosen with the Haar measure
  on $U(2N)$, and $U_2=V\Lambda V^*$, with $V$ a Haar distributed
  unitary independent of
  $U_{\infty}$, and $\Lambda$ a fixed symmetry of trace zero.

  As another example, specify to
  $a=
  \left[
    \begin{smallmatrix}
      0&t\\
      0&0
    \end{smallmatrix}
  \right].
  $
  As we will see, the spectrum and Brown measure are radially symmetric.
  The eigenvalues of $|\lambda-a|^2$ are
  \begin{equation}
    \label{eq:|lambda-nilpotent2|^2_Spectrum}
    \mu_\pm = \frac{t^2+2|\lambda|^2\pm t\sqrt{t^2+4|\lambda|^2}}%
                   {2}
  \end{equation}

  and hence, substituting this into \eqref{eq:a+u_SpectrumEquation}, we get
  $$
  \sigma(a+u) =
  \left\{
    \lambda
    :
    1-\frac{t^2}{2}
    \le |\lambda|^2
    \le \sqrt{\frac{t^2}{2}+\frac{1}{4}}+\frac{1}{2}
  \right\}
  $$
  which is a full disk for $t\ge \sqrt{2}$ and an annulus otherwise.
  For the density we substitute the parameters
  $$
  \|\lambda-a\|_2^2 = |\lambda|^2+\frac{t^2}{2}
  \qquad
  \det|\lambda-a|^2 = |\lambda|^4
  $$
  into formula \eqref{eq:a+u_Density}
  we get the radially symmetric density function
  $$
  p_{a+u}(\lambda)
  = \frac{4}{\pi}
    \left(
      \frac{ 2t^2}{ (4|\lambda|^2+t^2)^2}
      +
      \frac{(1-|\lambda|^2)^2 -  (1-2|\lambda|^2)t^2}
           { ((1-|\lambda|^2)^2 - t^2)^2 }
    \right)
  $$
\end{Example}

\subsection{An alternative expression for the Fuglede-Kadison determinant}
In order to treat more complicated examples,
instead of the integral \eqref{eq:Haar:logD:integral}
it will be more convenient to use
a more direct formula for the Kadison-Fuglede
determinant, which we state as a lemma.
\begin{Lemma}[{\cite[Proof of Theorem~4.4]{HaagerupLarsen:1999:BrownMeasure}}]
  Let $uh$ be an $R$-diagonal element and define functions on $\IR_+\setminus\{0\}$
  by
  \begin{align*}
    f(v) &= \tau((1+v h^2)^{-1})\\
    g(v) &= \frac{1-f(v)}{v f(v)}
  \end{align*}
  Then $g(v)$ is strictly decreasing with $g(]0,\infty[)=]\|h^{-1}\|_2^{-2},\|h\|_2^2[$
  and for every $z\in]\|h^{-1}\|_2^{-2},\|h\|_2^2[$ there is a unique $v>0$ such
  that $z^2=g(v)$. With this $v$ we have
  $$
  \log\Delta(uh-z)
  = \frac{1}{2}
    \int \log(1+v t)\,d\mu_{h^2}(t)
    +
    \frac{1}{2}
    \log\frac{z^2}{1+vz^2}
  $$
\end{Lemma}

For our problem of computing
$\log\Delta(\lambda-a-u) = \log\Delta(u^*(\lambda-a)-1)$
this translates as follows.
Putting $f(v,\lambda)= \tau((1+v|a-\lambda|^2)^{-1})$
and denoting $v(\lambda)$ the unique positive solution of the equation
$(1+v)f(v,\lambda)=1$, then
\begin{align*}
  \log\Delta(\lambda-a-u)
  &= \log\Delta(u^*(\lambda-a)-1) \\
  &= \frac{1}{2}\,\tau(\log(1+v|a-\lambda|^2)) - \frac{1}{2}\log(1+v)
  .
\end{align*}
Note that this approach cannot be used in the general setting of section~\ref{sec:a+uhBrown},
as it does not tell how to evaluate the Kadison-Fuglede determinant at $z=0$.

For the rest of this section we shall assume that $a$ is normal with spectral
measure $\mu_a$, so that we can write
\begin{equation}
  \label{eq:f(v,lambda)=int}
  f(v,\lambda)
  = \int \frac{d\mu_a(t)}{1+v|\lambda-t|^2}
\end{equation}
and again with $(1+v)f(v,\lambda)=1$,
\begin{align*}
  \log\Delta(a+u-\lambda)
  &=\frac{1}{2} \int\log(1+v|\lambda-t|^2)\, d\mu_a(t) - \frac{1}{2}\log(1+v)
\end{align*}
For the density of the Brown measure we obtain
\begin{align*}
  p(\lambda)
  &= \frac{2}{\pi}\,
     \partial_{\bar\lambda}\partial_\lambda\log\Delta(a+u-\lambda)\\
  &= \frac{1}{\pi}\,
     \partial_{\bar\lambda}
     \left(
       \int\frac{|\lambda-t|^2 \partial_\lambda v + v(\bar\lambda-\bar t)}%
                {1+v\,|\lambda-t|^2}
           \,d\mu(t)
       - \frac{1}{1+v}\,\partial_\lambda v
     \right)\\
  &= \frac{1}{\pi}\,
     \partial_{\bar\lambda}
     \left(
       \partial_\lambda v
       \underbrace{
         \left(
           \int \frac{|\lambda-t|^2}{1+v\,|\lambda-t|^2} \,d\mu(t) - \frac{1}{1+v}
         \right)
       }_{=0}
       + v \int \frac{\bar\lambda-\bar t}{1+v\,|\lambda-t|^2} \, d\mu(t)
     \right) \\
  &= \frac{1}{\pi}\,
     \partial_{\bar\lambda}
     \int
      \frac{1}{\lambda-t} \,
      \frac{v\,|\lambda-t|^2}{1+v\,|\lambda-t|^2} \, d\mu(t) \\
  &= \frac{1}{\pi}\,
     \partial_{\bar\lambda}
     \int
      \frac{1}{\lambda-t} \,
      \left(
        1 - \frac{1}{1+v\,|\lambda-t|^2} \, d\mu(t)
      \right) \\
  &= \frac{1}{\pi}
     \int
      \frac{1}{\lambda-t} \,
      \frac{|\lambda-t|^2 \,\partial_{\bar\lambda}v + v(\lambda-t)}%
           {(1+v\,|\lambda-t|^2)^2} \, d\mu(t) \\
  &= \frac{1}{\pi}\,
     \left(
       \partial_{\bar\lambda}v
       \int
        \frac{\bar\lambda-\bar t}
             {(1+v\,|\lambda-t|^2)^2} \, d\mu(t)
       +
       \int
        \frac{v}
             {(1+v\,|\lambda-t|^2)^2} \, d\mu(t)
  \right)
\end{align*}
now by implicit differentiation
\begin{align*}
  1 &= (1+v)f(v,\lambda) \\
  0 &= \partial_{\bar\lambda}v\,f(v,\lambda)
       +
       (1+v)(\partial_vf(v,\lambda)\,\partial_{\bar\lambda}v
             + \partial_{\bar\lambda}f(v,\lambda))\\
  \partial_{\bar\lambda}v
    &= - \frac{(1+v)\,\partial_{\bar\lambda}f}%
              {f+(1+v)\,\partial_vf} \\
  \partial_{\bar\lambda}f(v,\lambda)
    &= - \int \frac{v(\lambda-t)}{(1+v\,|\lambda-t|^2)^2}\,d\mu(t) \\
  \partial_{v}f(v,\lambda)
    &= - \int \frac{|\lambda-t|^2}{(1+v\,|\lambda-t|^2)^2}\,d\mu(t)
\end{align*}
and thus
\begin{equation}
  \label{eq:Haar:p(lambda)}
  p(\lambda)
   = \frac{1}{\pi}
     \left(
       \frac{1+v}{v\,(f(v,\lambda)+(1+v)\,\partial_vf(v,\lambda))}
       \,|\partial_\lambda f(v,\lambda)|^2
       +
       v f(v,\lambda)
       +v^2\,\partial_v f(v,\lambda)
  \right)
\end{equation}

\begin{figure}[pt]
    \begin{center} %
     \includegraphics[width=\textwidth, height=.44\textheight, keepaspectratio]{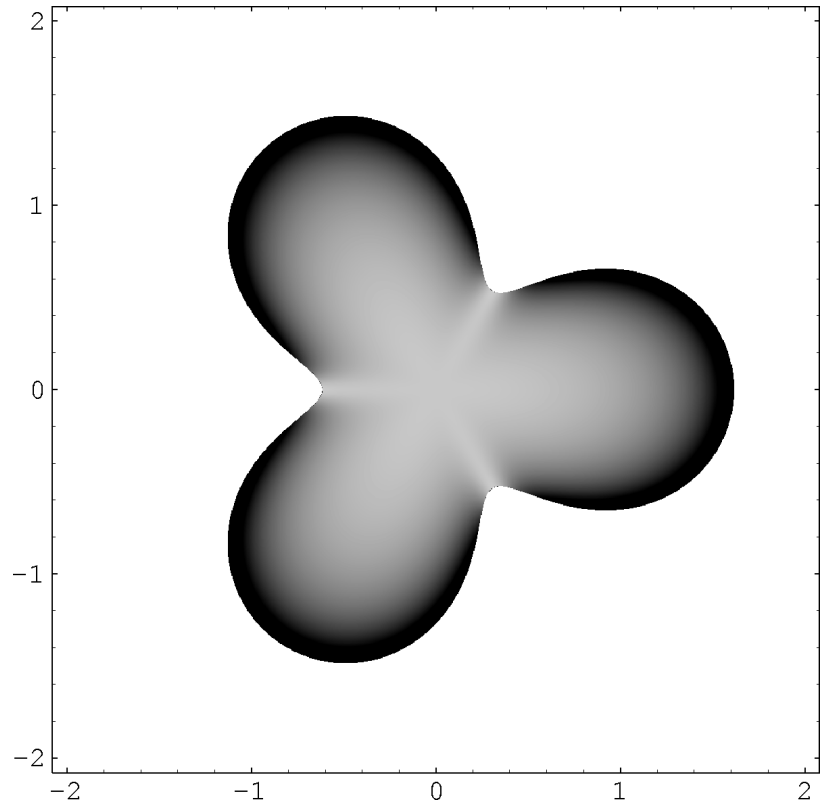}  
      \caption{Density of $\mu_{u_3+u_\infty}$}
      \label{fig:u3+uinfty:density}
    \end{center}
  \end{figure}

\begin{figure}[pb]
    \begin{center} %
     \includegraphics[width=\textwidth, height=.44\textheight, keepaspectratio]{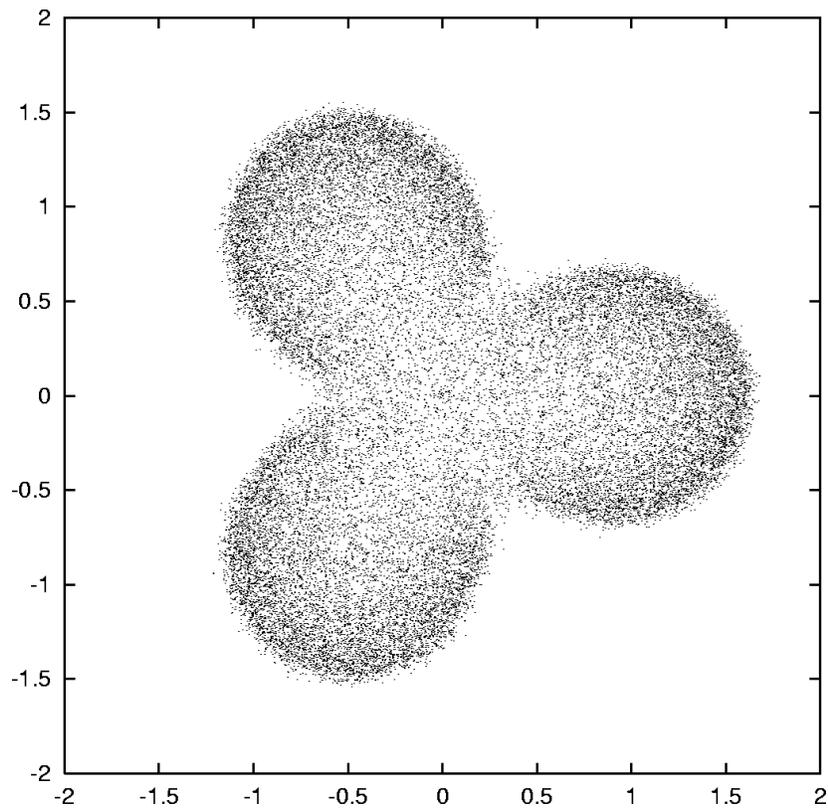}
      \caption{$200$ samples of eigenvalues of $150\times150$ random matrices $U_3+U_\infty$}
      \label{fig:U3+Uinfty_150}
    \end{center}
  \end{figure}

We will apply this in three situations here.
First consider a finite dimensional normal operator $a$, like e.g.\ $a=u_n$, the
generator of the von Neumann algebra of $\IZ_n$, then the integrals become
finite sums and can be evaluated numerically. As an example see fig.~\ref{fig:u3+uinfty:density},
which should again be compared to the corresponding samples of spectra of random matrices
in fig.~\ref{fig:U3+Uinfty_150}. There $U_3$ is a fixed $150\times150$ permutation matrix
with the same spectral distribution as $u_3$ and $U_\infty$ is again a $150\times 150$
standard unitary random matrix.

\begin{figure}[ht]
    \begin{center} %
     \includegraphics[width=\textwidth, height=.45\textheight, keepaspectratio]{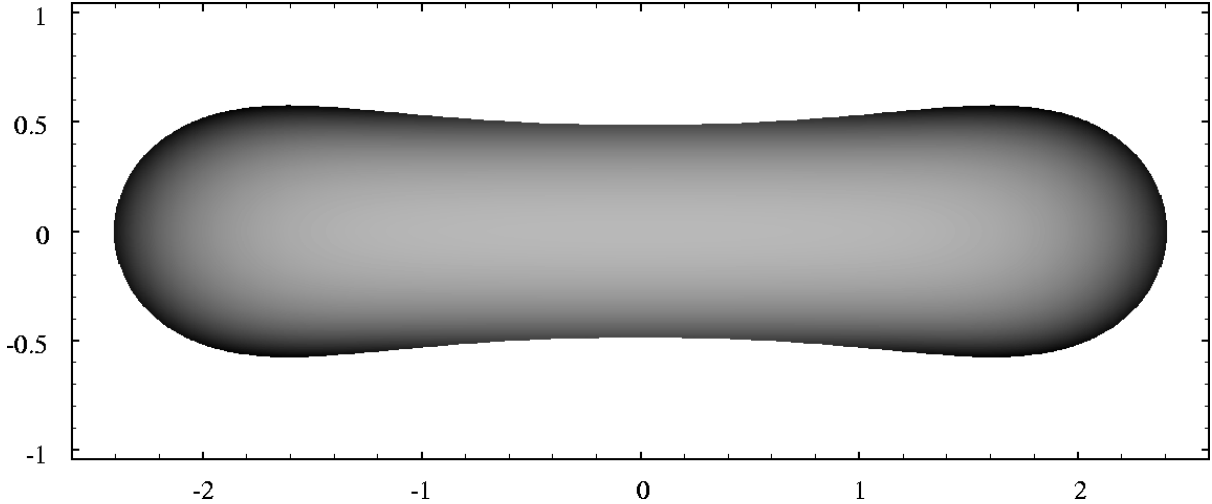}
      \caption{Density of $\mu_{u_2+v_2+u_\infty}$}
      \label{fig:u2+v2+uinfty:density}
    \end{center}
  \end{figure}

Secondly, assume that $a$ is self-adjoint.
Then we can factorize the denominator in the integral \eqref{eq:f(v,lambda)=int} as
$1+v\,|\lambda-t|^2= vt^2-v(\lambda+\bar\lambda)\,t+1+v|\lambda|^2 = v\,(t-z_0)(t-\bar z_0)$
where
$$
z_0 = \Re\lambda + \frac{i}{v}\sqrt{v^2(\Im\lambda)^2+v}
.
$$
From this we can express $f(v,\lambda)$ and therefore $p(\lambda)$ in terms of the
Cauchy transform $G(\zeta)$ of $a$ as follows.
\begin{align*}
  f(v,\lambda)
  &= \int\frac{d\mu(t)}{v\,(t-z_0)(t-\bar z_0)} \\
  &= \frac{1}{v}
     \int
      \frac{1}{z_0-\bar z_0}
      \left(
        \frac{1}{t-z_0}
        -
        \frac{1}{t-\bar z_0}
      \right)
      d\mu(t) \\
  &= -\frac{\Im G(z_0)}{\sqrt{v^2(\Im\lambda)^2+v}}
\end{align*}
As an example consider $a=u_2+v_2$, where $u_2$ and $v_2$ are the generators
of two free copies of $\IZ_2$. Then $a$ is self-adjoint and distributed according
to the arcsine law (or Kesten measure) and has Cauchy transform
$G(\zeta) = \frac{1}{\zeta\sqrt{1-\frac{1}{\zeta^2}}}$.
A picture of the density of the Brown measure of $u_2+v_2+u$ is presented in
fig.~\ref{fig:u2+v2+uinfty:density}.

\begin{figure}[ht]
    \begin{center} %
     \includegraphics[width=\textwidth, height=.45\textheight, keepaspectratio]{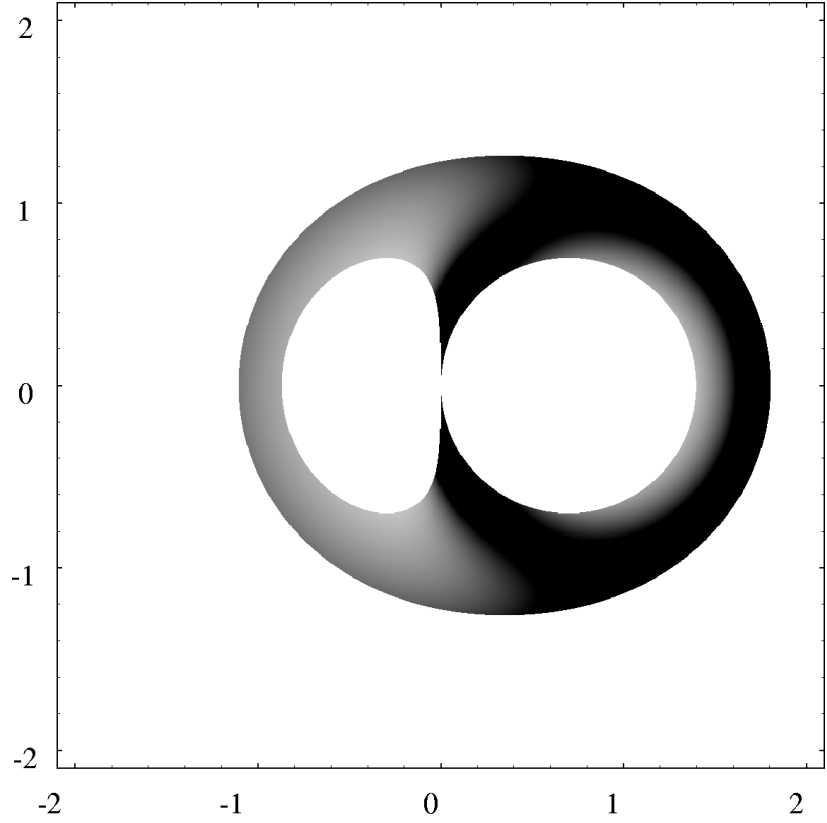}
      \caption{Density of $\mu_{u_q+u_\infty}$ at $q=0.7$}
      \label{fig:uq+uinfty:density}
    \end{center}
  \end{figure}

Finally, let us consider the free sum of an arbitrary unitary $v$ and a Haar unitary $u$.
Let $d\mu(\theta)$ be the spectral measure of $v$ on the unit circle.
For the evaluation of the integral \eqref{eq:f(v,lambda)=int} we factorize the
denominator again, this time writing
\begin{align*}
  f(v,\lambda)
  &= \int\frac{d\mu(\theta)}%
              {1+v\,|\lambda-e^{i\theta}|^2}\\
  &= \int
      \frac{d\mu(\theta)}%
           {1+v(|\lambda|^2+1)-v(\lambda e^{-i\theta}+\bar\lambda e^{i\theta})}\\
  &= -\int
      \frac{e^{i\theta}}%
           {v\bar\lambda e^{2i\theta} -(1+v(|\lambda|^2+1))e^{i\theta} +v\lambda}
      \,d\mu(\theta)\\
  &= -\frac{1}{v\bar\lambda}
     \int
      \frac{e^{i\theta}}%
           {(e^{i\theta}-z_+)(e^{i\theta}-z_-)}\\
\end{align*}
where
$$
z_\pm = \frac{1}{2v\bar\lambda}
        \left(
          1+v(|\lambda|^2+1)
          \pm
          \sqrt{(1+v(|\lambda|^2+1)^2)\, (1+v(|\lambda|^2-1)^2)}
        \right)
.
$$
Note that $|z_+ z_-|=|\frac{\lambda}{\bar\lambda}|$ and $|z_+|>|z_-|$, and thus
$|z_+|>1>|z_-|$.
\begin{align*}
  f(v,\lambda)
  &= \frac{1}{v\bar\lambda}
     \int
      \frac{e^{i\theta}}%
           {z_+z_-}
      \left(
        \frac{1}
        {(z_+-e^{i\theta})}
        -
        \frac{1}
        {(z_--e^{i\theta})}
      \right)
      d\mu(\theta) \\
  &= \frac{1}{v\bar\lambda}
     \int
      \frac{1}%
           {z_+z_-}
      \left(
        \frac{z_+}
        {(z_+-e^{i\theta})}
        -
        \frac{z_-}
        {(z_--e^{i\theta})}
      \right)
      d\mu(\theta) \\
  &= \frac{z_+G(z_+)-z_-G(z_-)}%
          {v\bar\lambda (z_+-z_-)} \\
  &= \frac{z_+G(z_+)-z_-G(z_-)}%
          {\sqrt{(1+v(|\lambda|^2+1)^2) \, (1+v(|\lambda|^2-1)^2)}}
\end{align*}
For the determination of the spectrum \eqref{eq:sigma(aplusu)} we need
\begin{align*}
  \|(\lambda-v)^{-1}\|_2^2
  &= \int \frac{d\mu(\theta)}{|\lambda-e^{i\theta}|^2} \\
  &= \frac{1}{|\lambda|^2-1}
     \int
      \left(
        \frac{\lambda}{\lambda-e^{i\theta}}
        +
        \frac{\bar\lambda}{\bar\lambda-e^{-i\theta}}
        -1
      \right)
      d\mu(\theta) \\
  &= \frac{\lambda\, G(\lambda)+\bar\lambda\, G(\bar\lambda)-1}%
          {|\lambda|^2-1}
\end{align*}
As an example let us consider for $q\in[-1,1]$ the unitary $u_q$ with Poisson distribution,
i.e.\ 
whose moments are $\tau(u_q^n)=q^{|n|}$. For $q=0$ this is the Haar distribution,
while for $q=1$ it is the Dirac measure at $1$. By Fourier transform, the density
of the spectral measure is
$$
d\mu_q(\theta) = \frac{1}{2\pi} \frac{1-q^2}{|1-q e^{i\theta}|^2}
.
$$
The Cauchy transform is 
$$
G_q(\zeta) = 
\begin{cases}
  \frac{1}{\zeta-q} & |\zeta|>1\\
  \frac{1}{\zeta-q^{-1}} & |\zeta|<1
\end{cases}
$$
and from this we get the other relevant functions
\begin{align*}
  \|(\lambda-u_q)^{-1}\|_2^2
  &= \begin{cases}
       \frac{|\lambda|^2-q^2}%
            {(|\lambda|^2-1)\,|\lambda-q|^2} & |\lambda| > 1\\
       \frac{q^{-2}-|\lambda|^2}%
            {(1-|\lambda|^2)\,|\lambda-q^{-1}|^2} & |\lambda| < 1
     \end{cases}\\
  f(v,\lambda)
  &= \frac{qz_+ - q^{-1}\,z_-}%
          {(z_+ - q)(z_- - q^{-1})}
     \frac{1}%
          {\sqrt{(1+v(|\lambda|^2+1)^2) \, (1+v(|\lambda|^2-1)^2)}}
\end{align*}
Substituting this into \eqref{eq:Haar:p(lambda)},
we get pictures like fig.~\ref{fig:uq+uinfty:density},
where $q=0.7$.

\section{Adding a circular element}

A standard circular element has the $*$-distribution of 
$C=S_1+iS_2$ where $S_1,S_2$ are
free standard 
semi-circular elements, i.e., self-adjoints whose distribution is
the semi-circle law ${1\over 2\pi}\sqrt{4-x^2}\,dx$ on $[-2,+2]$.
Its polar decomposition is $C=uh$ with $u$ a Haar unitary free with
$h$ (hence $C$ is 
$R$-diagonal), and $h$ has the quarter circular distribution 
${1\over \sqrt{2}\pi}\sqrt{8-x^2}\,dx$ on $[0,\sqrt 8]$. 
The symmetrized $\tilde h$ in Haagerup-Larsen's decomposition $C=a\tilde h$
has a semi-circular distribution of variance 2.
In this section we consider the Brown measure of $X_t=X_0+C_t$, where $X_0$ has
arbitrary $*$-distribution, it is free with $C_t$ and 
 $C_t$ is a circular element of variance $t$, i.e. $C_t\cong\sqrt{t\over 2}\, C$
where $C$ is a standard circular element. It will be convenient to assume that
the $C_t$ form a circular process, i.e., for each $s<t$, $C_t-C_s$ is $*$-free
with $C_s$. We shall use a heat equation like approach, by differentiating in
$t$.
One has
$$
\log\Delta(\lambda-X_t) = \frac{1}{2}\log\Delta(|\lambda-X_t|^2)
                  = \frac{1}{2}\lim_{\eps\to 0}\log\Delta(|\lambda-X_t|^2+\eps^2)
.
$$
Let us denote $H_{t,\eps}=|\lambda-X_t|^2+\eps^2$ and compute the derivative
$\frac{\partial}{\partial t}\log\Delta(H_{t,\eps})$.
To this end let $dt$ be small, $dC_t=C_{t+dt}-C_t$
(so that $\tau(dC_t^*dC_t)=dt$).
Then
\begin{align*}
  H_{t+dt,\eps}
  &= |\lambda-X_{t+dt}|^2+\eps^2 \\
  &= |\lambda-X_t-dC_t|^2 + \eps^2\\
  &= |\lambda-X_t|^2 - (\lambda-X_t)^*dC_t-dC_t^*(\lambda-X_t)+|dC_t|^2+\eps^2\\
  &= H_{t,\eps} - (\lambda-X_t)^*dC_t-dC_t^*(\lambda-X_t)+dC_t^*dC_t\\
  &= H_{t,\eps}\left[
  1- H_{t,\eps}^{-1}( (\lambda-X_t)^*dC_t+dC_t^*(\lambda-X_t)-dC_t^*dC_t)
  \right]\\
\intertext{and hence}
  \log\Delta(H_{t+dt,\eps})
  &= \log\Delta(H_{t,\eps})
     +
     \log\Delta(1- H_{t,\eps}^{-1}( (\lambda-X_t)^*dC_t+dC_t^*(\lambda-X_t)-dC_t^*dC_t))\\
  &= \log\Delta(H_{t,\eps})
     +
     \tau(\log|1- H_{t,\eps}^{-1}( (\lambda-X_t)^*dC_t+dC_t^*(\lambda-X_t)-dC_t^*dC_t)|)
  .
\end{align*}
Now observe that 
\begin{align*}
  \tau(\log|1+a\sqrt{dt} + b\,dt|)
  &= \frac{1}{2}
     \tau(\log|1+a\sqrt{dt} + b\,dt|^2)\\
  &= \frac{1}{2}
     \tau\left(
           \log(1+(a+a^*)\sqrt{dt} + (b+b^*)dt + a^*a\,dt + \Ord((dt)^{3/2}) 
     \right)\\
  &= \frac{1}{2}
     \tau
     \left(
       (a+a^*)\sqrt{dt}
       +
       (b+b^*+a^*a)dt
       -
       \frac{1}{2}
       (a+a^*)^2 dt
     \right)
     +
     \Ord((dt)^{3/2})\\
  &= \frac{1}{2}
     \tau
     \left(
       (a+a^*)\sqrt{dt}
       +
       (b+b^*-\frac{a^2+{a^*}^2}{2})dt
     \right)
     +
     \Ord((dt)^{3/2})
\end{align*}
In our situation we have
$
a=-H_{t,\eps}^{-1}
  \left(
    (\lambda-X_t)^*\frac{dC_t}{\sqrt{dt}}
    +
    \frac{dC_t}{\sqrt{dt}}^*(\lambda-X_t)
  \right)
$
and $b=H_{t,\eps}^{-1}\frac{dC_t^*dC_t}{dt}$, 
so that $\tau(a)=0$ and $\tau(b)=\tau(H_{t,\eps}^{-1})$ by freeness of $dC_t$ and
$\{H_{t,\eps}, \lambda-X_t\}$. Further we have
\begin{align*}
  \tau(a^2)
  &= \tau\left(
           \left(
             H_{t,\eps}^{-1} (\lambda-X_t)^*\, \frac{dC_t}{\sqrt{dt}}
           \right)^2
           +
           \left(
             H_{t,\eps}^{-1}\, \frac{dC_t^*}{\sqrt{dt}}\, (\lambda-X_t)
           \right)^2
\right.\\&\hskip3em\left.\phantom{\left(\frac{dC_t^*}{\sqrt{dt}}\right)^2}
           +
           2
           H_{t,\eps}^{-1} 
           (\lambda-X_t)^*
           \,
           \frac{dC_t}{\sqrt{dt}}
           \,
           H_{t,\eps}^{-1}
           \,
           \frac{dC_t^*}{\sqrt{dt}}
           \,
           (\lambda-X_t)
         \right)
\intertext{
  and using the formula
  $\tau(a_1b_1a_2b_2)=\tau(a_1)\tau(a_2)\tau(b_1b_2)$
  if $\{a_1,a_2\}$ is free from $\{b_1,b_2\}$ and $\tau(b_1)=\tau(b_2)=0$,
  we see that only the last term is nonzero and equal to
}
  &= 2\tau\left(
            (\lambda-X_t)
            \,
            H_{t,\eps}^{-1}
            \,
            (\lambda-X_t)^*
            \,
            \frac{dC_t}{\sqrt{dt}}
            \,
            H_{t,\eps}^{-1}
            \,
            \frac{dC_t^*}{\sqrt{dt}}
          \right) \\
  &= 2
     \tau( (\lambda-X_t)
           \,
           H_{t,\eps}^{-1}
           \,
           (\lambda-X_t)^* )
     \,
     \tau( H_{t,\eps}^{-1} ) \\
  &= 2\,\tau(H_{t,\eps}^{-1}\, (H_{t,\eps}-\eps^2)\,
     \tau( H_{t,\eps}^{-1} ) \\
  &= 2\,\tau(H_{t,\eps}^{-1})
     +
     2\eps^2 \tau(H_{t,\eps}^{-1})^2
\end{align*}
so that 
\begin{align*}
  \frac{\log\Delta(H_{t+dt,\eps}) -\log\Delta(H_{t,\eps}) }{dt}
  &= \frac{1}{2}
     \,
     (2\tau(H_{t,\eps}^{-1})
     - 
     2\tau(H_{t,\eps}^{-1})
     +
     2\eps^2 \tau(H_{t,\eps}^{-1})^2
     +
     \Ord((dt)^{1/2})\\
  &= \eps^2\tau(H_{t,\eps}^{-1})^2 + \Ord((dt)^{1/2})
\end{align*}

hence
$$
\frac{\partial}{\partial t}
\log\Delta(H_{t,\eps}) = \eps^2\tau(H_{t,\eps}^{-1})^2
$$
and
$$
\log\Delta(H_{t,\eps}^{-1})
= \log\Delta(H_{0,\eps}^{-1})
  +
  \int_0^t
   \eps^2
   \,
   \tau(H_{s,\eps}^{-1})^2
   \,
   ds
.
$$
Let $a_{\lambda,s}$ be a  self-adjoint element with symmetric distribution,
 whose absolute value is distributed
as $\vert\lambda-X_s\vert$.
Now note that by the Stieltjes inversion formula
\begin{align*}
  \eps\,\tau(H_{s,\eps}^{-1})
  &= \tau\left(
           \eps\,
           (|\lambda-X_s|^2+\eps^2)^{-1}
         \right)\\
  &= -\tau\left(
           \Im[(i\eps-a_{\lambda,s})^{-1}]
         \right)\\
  &\rightarrow_{\eps\to0}
     \pi\,
     \left.
       \frac{d\mu_{a_{\lambda,s}}(x)}%
            {dx}
     \right|_{x=0}
\end{align*}
i.e., the density at $0$ of the distribution of
$a_{\lambda,s}$.
Now we need the following
\begin{Lemma}
Let $a$ be a self adjoint symmetrically distributed element, free with $S$ and
$C$, where $S$ and $C$ are a semi-circular and a circular element of same
variance
respectively,
then $\vert a+S\vert$ and $\vert \vert a\vert+C\vert$ 
have the same distribution.
\end{Lemma}
\begin{proof}
  Let $b$ be a symmetry free with $\{a,S,C\}$, then by
  {\cite[Prop.~4.2]{HaagerupLarsen:1999:BrownMeasure}}
  $b a$ and $b S$ are $*$-free, thus 
  $\vert a+S\vert=\vert b a+b S\vert$ is distributed as $\vert b a+C\vert$. 
  Now using the fact that multiplying with a free Haar unitary $u$ 
  does not change
  the $*$-distribution of $C$, we can replace the latter according to 
  $C\cong u^* C$,
  and get the following equalities of $*$-distributions 
  $$
  |b a+C|\cong |b a+u^*C|
  \cong |ub a+C|
  \cong \bigl|
  u|a|+C
  \bigr|
  \cong \bigl|
  |a|+C
  \bigr|
  $$
\end{proof}

Using the lemma we get
$$
|\lambda-X_s|=\vert\lambda-X_0-C_s\vert \cong \bigl|
                      a_{\lambda} + S_s
                    \bigr|
$$
where $a_{\lambda}$ is the symmetrization of $\vert\lambda-X_0\vert$, free with 
the semicircular $S_s$ and therefore
$$
|a_{\lambda,s}| \cong |a_{\lambda}+S_s|
$$
It follows from Corollary 3 of
\cite[p.~711]{Biane:1997:SemiCircularConvolution} that the distribution of
$a_{\lambda,s}$
has a
 density at $0$ which  is
$p_s(0) = \frac{v(s)}{\pi s}$,
with
\begin{equation}
  \label{eq:SemiCirclev(s)}
  v(s) = \inf\left\{
              v\ge 0 : \int
              \frac{d\mu_{\vert\lambda-X_0\vert}(x)}{x^2+v^2}\le\frac{1}{s}
            \right\}
\end{equation}
If $\lambda\not \in \sigma(X_0)$, then by e.g. 
\cite{Biane:1997:SemiCircularConvolution} 
$$
  \begin{aligned}
     \eps\tau(H_{s,\eps}^{-1})
      &\leq \sup_{x\in \IR}{d\mu_{a_{\lambda,s}}(x)\over dx}\\
      &\leq {1\over \pi\sqrt s}
 \end{aligned}     
$$
furthermore for $s$ small enough, $\lambda\not\in\sigma(X_s)$ and
$\tau(|\lambda-X_s\lambda|^{-2})$ is bounded above, hence
$\eps\tau(H_{s,\eps}^{-1})$ also, 
therefore we can apply the dominated convergence theorem 
and  we get 
\begin{equation}
  \label{eq:logDeltaXt}
  \begin{aligned}
    \log\Delta(\lambda-X_t)
    &= \frac{1}{2}
       \lim_{\eps\to 0}
       \log\Delta(H_{0,\eps})
       +
       \frac{1}{2}\int_0^t
        \eps^2\tau(H_{s,\eps}^{-1})^2 ds\\
    &= \log\Delta(\lambda-X_0)
       +
       \frac{1}{2}
       \int_0^t
        \frac{v(s)^2}{s^2}
        \,ds\\
    &= \log\Delta(\lambda-X_0)
       +
       \frac{1}{2}
       \int_{t_{\lambda}}^t
        \frac{v(s)^2}{s^2}
        \,ds\\
  \end{aligned}
\end{equation}
where $t_\lambda = \inf\{t : v(t)>0\} =
 (\int\frac{d\mu_{\vert\lambda-X_0\vert}(x)}{x^2})^{-1}$.
So whenever $\lambda\not\in\sigma(X_0)$, 
the density of the Brown measure is
\begin{align*}
  p_{\lambda-X_t}(\lambda)
  &= \frac{1}{\pi}
     \partial_{\bar{\lambda}}\partial_\lambda
     \int_{t_\lambda}^t 
      \frac{v(s)^2}{s^2}
      \,ds \\
  &= \frac{1}{\pi}
     \,
     \partial_{\bar{\lambda}}
     \left(
       \int_{t_\lambda}^t
        \frac{\partial_{\lambda}v(s)^2}{s^2}
        \,ds
       -
       \frac{v(t_\lambda)^2}%
            {t_\lambda^2}
       \partial_\lambda t_\lambda
     \right)
\end{align*}
and the second summand will be zero if $v(t)$ is continuous at $t_\lambda$.

\begin{Example}[$2\times2$ matrix]
  Let $X_0=a$ be as in example 4.1, and consider $X_t=a+C_t$.
  Let again $\mu_\pm$ be the eigenvalues of $(\lambda-a)^*(\lambda-a)$,
  then the relevant parameters are
  \begin{align*}
    \|\lambda-a\|_2^2
    &= \frac{\mu_+ + \mu_-}{2}\\
    \|(\lambda-a)^{-1}\|_2^2
    &= \frac{1}{2}
        \left(
          \frac{1}{\mu_+}
          +
          \frac{1}{\mu_-}
        \right)
     = \frac{\|\lambda-a\|^2_2}{\det|\lambda-a|^2} \\
    t_\lambda
    &= \left(
        \int
        \frac{d\mu_{|\lambda-a|}(x)}%
             {x^2}
        \,dx
      \right)^{-1}
    = \|(\lambda-a)^{-1}\|_2^{-2}
    =  \frac{\det|\lambda-a|^2}{\|\lambda-a\|^2_2}
    .
  \end{align*}
  The function $v(s)^2$ is the solution of the quadratic equation
  \begin{align*}
    \frac{1}{s}
    &= \frac{1}{2}
       \left(
         \frac{1}{\mu_++v^2}
         +
         \frac{1}{\mu_-+v^2}
       \right) \\
    &= \frac{1}{2}
       \frac{\mu_++\mu_-+2v^2}%
            {\mu_+\mu_-+(\mu_++\mu_-)v^2+v^4}\\
    &= \frac{\|\lambda-a\|^2_2+v^2}%
            {\det|\lambda-a|^2 + 2\|\lambda-a\|_2^2\, v^2 + v^4}
  \end{align*}
  which is explicitly
  \begin{align*}
  v(s)^2
  &= \frac{1}{2}
     \left(
       s - 2\|\lambda-a\|_2^2
       \pm \sqrt{(s- 2\|\lambda-a\|_2^2)^2
                  - 4(\det|\lambda-a|^2 - s\|\lambda-a\|_2^2 )}
     \right)\\
  &= \frac{1}{2}
     \left(
       s - 2\|\lambda-a\|_2^2 \pm \sqrt{s^2 + 4(\|\lambda-a\|_2^4-\det|\lambda-a|^2)}
     \right)
  \end{align*}
  Now we have to choose the right branch of the square root.
  To this end, let us compute the spectrum of $X_t$:
  Assume $\lambda\not\in\sigma(a)$, then
  $\lambda\in\sigma(a+C_t)$ if and only if 
  $1-C_t(\lambda-a)^{-1}$ is not invertible.
  Now $C_t(\lambda-a)^{-1}$ is $R$-diagonal and not invertible, so by
  Theorem~\ref{thm:HaagerupLarsen:5.1}~\eqref{item:HaagerupLarsen:5.1:notinvertible},
  $1$ is in its spectrum if and only if its spectral radius is at least $1$
  and using Proposition~\ref{prop:HaagerupLarsen:4.10} we get the inequality
  $$
  1 \le \rho(C_t(\lambda-a)^{-1}) = \|C_t\|_2 \, \|(\lambda-a)^{-1}\|_2
  $$
  in other words,
  $$
  \det|\lambda-a|^2 \le t \, \|\lambda-a\|_2^2
  $$
  and hence for $s<t$, $\det|\lambda-a|^2-s\,\|\lambda-a\|_2^2 < 0$,
  only the ``$+$'' branch gives a nonnegative solution.
  Consequently
  \begin{align*}
    \log\Delta(\lambda-X_t)&-\log\Delta(\lambda-X_0) \\
    &= \frac{1}{2}
       \int_{t_\lambda}^t
        \frac{1}{2s}
        -
        \frac{\|\lambda-a\|_2^2}{s^2}
        +
        \frac{\sqrt{s^2+4(\|\lambda-a\|_2^4-\det|\lambda-a|^2)}}%
             {2s^2}
        \,ds \\
    &= \left.
         \frac{1}{4}
         \log s
         +
         \frac{\|\lambda-a\|_2^2}{2s}
         +
         \frac{1}{4}
         \log\left(
               s+\sqrt{s^2+4(\|\lambda-a\|_2^4-\det|\lambda-a|^2)}
             \right)
\right.\\ &\hskip5em\left.
         -
         \frac{1}{4s}
         \sqrt{s^2+4(\|\lambda-a\|_2^4-\det|\lambda-a|^2)}
       \right|_{s=t_\lambda}^t
  \end{align*}
Now observe that
$$
\sqrt{t_\lambda^2+4(\|\lambda-a\|_2^4-\det|\lambda-a|^2)}
= \frac{2\|\lambda-a\|_2^4 - \det|\lambda-a|^2}%
            {\|\lambda-a\|_2^2}\\
$$
and hence, denoting
\begin{equation}
  \label{eq:2x2+circular:R(lambda)}
  R(\lambda) = 4(\|\lambda-a\|_2^4-\det|\lambda-a|^2) = (\mu_+-\mu_-)^2  
\end{equation}
we get
\begin{align*}
  \log&\Delta(\lambda-X_t) - \log\Delta(\lambda-X_0)\\
      &= \frac{1}{4}
         \log t
         +
         \frac{\|\lambda-a\|_2^2}{2t}
         +
         \frac{1}{4}
         \log\left(
               t+\sqrt{t^2+R(\lambda)}
             \right)
         -
         \frac{1}{4t}
         \sqrt{t^2+R(\lambda)}
\\ &\phantom{=+}
         -
         \frac{1}{4}
         \log \frac{\det|\lambda-a|^2}{\|\lambda-a\|_2^2}
         -
         \frac{\|\lambda-a\|_2^4}{2\det|\lambda-a|^2}
         -
         \frac{1}{4}
         \log\frac{\det|\lambda-a|^2+2\|\lambda-a\|_2^4-\det|\lambda-a|^2}%
                  {\|\lambda-a\|_2^2}
\\ &\phantom{=+}
         +
         \frac{\|\lambda-a\|_2^2}%
              {4\det|\lambda-a|^2}
         \,
         \frac{2\|\lambda-a\|_2^4-\det|\lambda-a|^2}%
              {\|\lambda-a\|_2^2} \\
      &= \frac{1}{4}
         \log t
         +
         \frac{\|\lambda-a\|_2^2}{2t}
         +
         \frac{1}{4}
         \log\left(
               t+\sqrt{t^2+R(\lambda)}
             \right)
         -
         \frac{1}{4t}
         \sqrt{t^2+R(\lambda)}
\\ &\phantom{=+}
         -
         \frac{1}{4}
         \log\det|\lambda-a|^2 
         -
         \frac{1}{4}
         \log 2
         -
         \frac{1}{4}
\end{align*}
and finally the density is
(note that 
$\partial_{\bar{\lambda}} \partial_\lambda \log\det|\lambda-a|^2=0$
and $\partial_{\bar{\lambda}} \partial_\lambda \|\lambda-a\|_2^2=1$)
\begin{align*}
  p_{a+C_t}(\lambda)  
  &= \frac{2}{\pi}\,
     \partial_{\bar{\lambda}} \partial_\lambda
     \log\Delta(\lambda-X_t) \\
  &= \frac{1}{\pi t}
     +
     \frac{1}{2\pi}\,
     \partial_{\bar{\lambda}} \partial_\lambda
     \left(
       \log\left(
             t+\sqrt{t^2+R(\lambda)}
           \right)
       -
       \frac{\sqrt{t^2+R(\lambda)}}{t}
     \right) \\
  &= \frac{1}{\pi t}
     +
     \frac{1}{2\pi}\,
     \partial_{\bar{\lambda}}
     \left(
       \frac{1}%
            {t+\sqrt{t^2+R(\lambda)}}
       -
       \frac{1}{t}
     \right)
     \frac{\partial_\lambda R(\lambda)}%
          {2\sqrt{t^2+R(\lambda)}} \\
  &= \frac{1}{\pi t}
     +
     \frac{1}{4\pi}\,
     \partial_{\bar{\lambda}}
     \left(
       \frac{t-\sqrt{t^2+R(\lambda)}}
            {tR(\lambda)}
       \partial_\lambda R(\lambda)
     \right) \\
  &= \frac{1}{\pi t}
     +
     \frac{1}{4\pi t}
     \left(
       -
       \frac{|\partial_\lambda R(\lambda)|^2}%
            {2R(\lambda) \sqrt{t^2+R(\lambda)}}
       +
       \left(
         t-\sqrt{t^2+R(\lambda)}         
       \right)
       \frac{R(\lambda)\partial_{\bar{\lambda}}\partial_\lambda R(\lambda)
               -|\partial_\lambda R(\lambda)|^2 }%
            {R(\lambda)^2}
     \right)
\end{align*}
Again we can specify to $
a=
\left[
  \begin{smallmatrix}
    0&1\\
    1&0
  \end{smallmatrix}
\right]
$,
and we get the spectrum 
$$
\sigma(a+C_t) = \{\lambda : |\lambda^2-1|^2 \le t (|\lambda|^2+1) \}
$$
Note that for $t=1$ this is the same as $\sigma(u_2+u)$ from example 4.1.
However this time the density is a function of the real part alone,
namely substituting $\mu_\pm = |\lambda\pm1|^2$ into \eqref{eq:2x2+circular:R(lambda)},
we get $R(\lambda)=4(\lambda+\bar\lambda)^2$
and consequently the density depends only on the real part
$$
p_{a+C_t}(x+iy)
=\frac{1}{\pi t}
 +
 \frac{1}{8\pi x^2}
 \left(
   \frac{t}{\sqrt{t^2+16x^2}}-1
 \right)
$$

The situation for the nilpotent $2\times2$ matrix
$a=
\left[
  \begin{smallmatrix}
    0&1\\
    0&0
  \end{smallmatrix}
\right]
$
is as follows.
We have computed the eigenvalues of $|\lambda-a|^2$ in \eqref{eq:|lambda-nilpotent2|^2_Spectrum},
and thus
$$
\sigma(a+C_t)
= \left\{
    \lambda
    :
    2|\lambda|^4
    \le
    t\, (1+2|\lambda|^2)
  \right\}
$$
which is the disk with radius $\sqrt{\sqrt{\frac{t^2}{4}+\frac{1}{2}}+\frac{t}{2}}$.
This is the same as $\sigma(a+ \sqrt{t}\,u)$, but with the possible hole removed.
Furthermore we get $R(\lambda)=(\mu_+-\mu_-)^2=1+4|\lambda|^2$ and the density function is
again rotationally symmetric:
$$
p_{a+C_t}(\lambda)
= \frac{1}{\pi t}
  \left(
    1
    -
    \frac{2|\lambda|^2}%
         {(1+4|\lambda|^2)\sqrt{t^2+1+4|\lambda|^2}}
    +
    \frac{t-\sqrt{t^2+R}}%
         {(1+4|\lambda|^2)^2}
  \right)
$$
\end{Example}

\begin{Example}[Elliptic law]
  An interesting example is given by the so-called \emph{elliptic random variable}
  $S_{\alpha}+i S_{\beta}$,
  where $S_{\alpha}$ and $S_{\beta}$ are free semicircular variables of
  variances $\alpha$ and $\beta$. 
  Note that for $\alpha=\beta$ this is a circular variable $C_{2\alpha}$.
  The Brown measure has been computed by Haagerup (unpublished) by another
  method.
  The name \emph{elliptic} stems from the shape of its spectrum, 
  which is an ellipse.
  This can be seen as follows.
  Assuming that $\alpha>\beta$ let $\gamma=\alpha-\beta$,
  then for $\lambda\not\in \sigma(S_\gamma)=[-2\sqrt{\gamma},2\sqrt{\gamma}]$ we have
  $\lambda\in\sigma(S_\gamma+C_{2\beta})$ if and only if
  $1-C_{2\beta}(\lambda-S_\gamma)^{-1}$ is not invertible.
  From Theorem~\ref{thm:HaagerupLarsen:5.1} we infer that the spectrum of 
  $C_{2\beta}(\lambda-S_\gamma)^{-1}$ is the disk centered at zero with radius
  $\|C_{2\beta}(\lambda-S_\gamma)^{-1}\|_2$, so that we get
  $$
  \sigma(S_\gamma+C_{2\beta})
  = \{ \lambda : 1 \le 2\beta\, \|(\lambda-S_\gamma)^{-1}\|_2^2 \}
  $$
  We use formula~\eqref{eq:||(l-a)^(-1)||_2} for the Cauchy transform
  $G_{S_\gamma}(\zeta) = \frac{\zeta-\sqrt{\zeta^2-4\gamma}}{2\gamma}$ to get
  \begin{equation}
    \label{eq:||l-S^(-1)||}
  \|(\lambda-S_\gamma)^{-1}\|_2^2
  = \frac{1}{2\gamma}
    \left(
      \frac{\sqrt{\lambda^2-4\gamma} - \sqrt{\bar{\lambda}^2-4\gamma}}%
           {\lambda-\bar{\lambda}}
      -
      1
    \right)
  \end{equation}
  and hence the spectrum is
  $$
  \left\{
    \lambda
    : 
    \frac{\sqrt{\lambda^2-4\gamma} - \sqrt{\bar{\lambda}^2-4\gamma}}%
         {\lambda-\bar{\lambda}}
    \ge \frac{\gamma+\beta}{\beta} = \frac{\alpha}{\beta}
  \right\}
  $$
  Now consider the Zhukowski transformation $f:\xi\mapsto\frac{1}{\xi}+\gamma\xi$,
  which maps the circles
  $\{\frac{e^{i\theta}}{t} : 0\le\theta<2\pi\}$ to the ellipses
  $
  \left\{
    \left(
      \frac{\gamma}{t}+t
    \right)
    \cos\theta
    +
    i
    \left(
      \frac{\gamma}{t}-t
    \right)
    \sin\theta
    :
    0\le \theta<2\pi
  \right\}
  $
  and hence the open disk
  $\{\xi:|\xi|<\frac{1}{\sqrt{\gamma}}\}$
  bijectively onto $\IC\setminus[-2\sqrt{\gamma},2\sqrt{\gamma}]$.
  Note that the excluded interval is exactly the spectrum of $S_\gamma$.
  So assume that $\lambda=f(\xi)$ with $|\xi|<\frac{1}{\sqrt{\gamma}}$
  is not in the spectrum of $S_\gamma$, then
  observe that
  $$
  \lambda^2-4\gamma
  = \frac{1}{\xi^2}+2\gamma+\gamma^2\xi^2 - 4\gamma
  = \left(
      \frac{1}{\xi} - \gamma\xi
    \right)^2
  $$
  and hence $\lambda\in\sigma(S_\gamma+C_{2\beta})$ if and only if
  $$
  \frac{\alpha}{\beta}
  \le \frac{\frac{1}{\xi}-\frac{1}{\bar{\xi}} -\gamma\xi +\gamma\bar{\xi} }%
           {\frac{1}{\xi}-\frac{1}{\bar{\xi}} +\gamma\xi -\gamma\bar{\xi} }
  = \frac{1+\gamma|\xi|^2}{1-\gamma|\xi|^2}
  .
  $$
  This inequality reduces to
  $$
  |\xi|^2\ge \frac{1}{\alpha+\beta}
  ,
  $$
  thus 
  $$
  \sigma(S_\gamma+S_{2\beta})\setminus[-2\sqrt{\gamma},2\sqrt{\gamma}]
  = \left\{
      f(\xi) : \frac{1}{\sqrt{\alpha+\beta}} \le |\xi| < \frac{1}{\sqrt{\gamma}}
    \right\}
  $$
  and taking the closure of this set we obtain $\sigma(S_\gamma+C_{2\beta})$
  as the interior of the ellipse
  \begin{equation}
    \label{eq:sigma(Salpha+iSbeta)}
    \left\{
      \frac{2\alpha}{\sqrt{\alpha+\beta}}
      \cos\theta
      +
      \frac{2\beta}{\sqrt{\alpha+\beta}}
      \,
      i\sin\theta
      :
      0\le\theta<2\pi
    \right\}
    .
  \end{equation}
  Now let us turn to the Brown measure.
  As already noted, the method from section~\ref{sec:a+uhBrown} will not work
  on $a=S_\gamma$.
  Indeed the $R$-transform of $|\lambda-S_\gamma|^2$ can be computed from the inverse of
  $$
  G_{|\lambda-S_\gamma|^2}(\zeta)
  = \frac{1}{2\gamma}
    \left(
      1
      -
      \frac{\sqrt{x_+^2-4\gamma} - \sqrt{x_-^2-4\gamma}}%
           {x_+ - x_-}
    \right)
  $$
  where $x_\pm$ are as in \eqref{eq:Gl-a2:xplusminus}.
  Let $\lambda = \xi+ i \eta$, then
  we can rewrite $x_\pm = \xi\pm\sqrt{\zeta-\eta^2}$
  and abbreviating $y=\sqrt{\zeta-\eta^2}$,
  solve the equation $G_{|\lambda-S_\gamma|^2}(\zeta) = z$ for $y$, which gives
  $$
  y^2 = \frac{\xi^2}{(1-2\gamma z)^2}
        +
        \frac{1}{z(1-\gamma z)}
  .
  $$
  It follows that
  $K(z) = y^2+\eta^2$
  and
  $$
  R_{|\lambda-S_\gamma|^2}(z)
       = z\,K(z) - 1
       = \frac{\gamma z}{1-\gamma z}
         +
         \frac{\xi^2 z}{(1-2\gamma z)^2}
         +
         \eta^2 z
   ;
  $$
  for real $\lambda$ this has been used in
  \cite{HiwatashiNagisaKurodaYoshida:1999:chisquare} to characterize
  the semicircular distributions.
  In order to get the determining series $f_{u|\lambda-S_\gamma|}$
  according to \eqref{eq:fxl-1=Rxl2} one has to solve a fourth order equation,
  which is not suitable for further computations.
  So we have to use formula~\eqref{eq:logDeltaXt},
  for which we need $v(s)$ from \eqref{eq:SemiCirclev(s)} first.
  We have done most of the work already, since
  $\int\frac{d\mu(x)}{|\lambda-x|^2+v^2} = - G_{|\lambda-S_\gamma|^2}(-v^2)$,
  thus
  $$
  v(s)^2 = - K_{|\lambda-S_\gamma|^2}(-\frac{1}{s})
       = - \left(
             \frac{\xi^2 s^2}{(s+2\gamma)^2}
             -
             \frac{s^2}{s+\gamma}
             +
             \eta^2
           \right)
  $$
  and
  $$
  \frac{v(s)^2}{s^2}
  = - \frac{(\lambda+\bar\lambda)^2}{4(s+2\gamma)^2}
    + \frac{1}{s+\gamma}
    -
    \frac{(\lambda-\bar\lambda)^2}{4s^2}
  $$
  and the density becomes, with
  \begin{align}
    p_{S_\alpha+iS_\beta}(\lambda)
    &= \frac{1}{\pi}
       \,
       \partial_{\bar\lambda}
       \int_{t_\lambda}^{2\beta}
        \frac{\partial_\lambda v(s)^2}{s^2}
        \,ds \nonumber\\
    &= \frac{1}{\pi}
       \,
       \partial_{\bar\lambda}
       \int_{t_\lambda}^{2\beta}
        \left(
          -\frac{2(\lambda+\bar\lambda)}%
                {4(s+2\gamma)^2}
          +\frac{2(\lambda-\bar\lambda)}{4s^2}
        \right)
        \,ds     \nonumber\\
    &= \frac{1}{2\pi}
       \,
       \partial_{\bar\lambda}
       \left.
         \left(
           \frac{\lambda+\bar\lambda}{s+2\gamma}
           -
           \frac{\lambda-\bar\lambda}{s}
         \right)
       \right|_{t_{\lambda}}^{2\beta}\nonumber\\
    &= \label{eq:pS+iSdlb}
       \frac{1}{4\pi}
       \left(
         \frac{1}{\alpha}
         +
         \frac{1}{\beta}
       \right)
       -\frac{1}{2\pi}
       \partial_{\bar\lambda}
         \left(
           \frac{\lambda+\bar\lambda}{t_\lambda+2\gamma}
           -
           \frac{\lambda-\bar\lambda}{t_\lambda}
         \right)
  \end{align}
  Now $t_\lambda = \|(\lambda-S_\gamma)^{-1}\|_2^{-2}$ has been computed above in
  \eqref{eq:||l-S^(-1)||},
  and denoting $\omega = \sqrt{\lambda^2-4\gamma}$, it is
  $$
  t_\lambda = 2\gamma
              \left(
                \frac{\omega-\bar\omega}{\lambda-\bar\lambda}
                -1
              \right)^{-1}
  $$
  and we claim now that the second summand in \eqref{eq:pS+iSdlb}
  is zero. For this note that
  $\omega^2-{\bar\omega}^2 = \lambda^2 - \bar\lambda^2$ and hence
  \begin{align*}
    -
    \frac{1}{2\pi}
    \,
    \partial_{\bar\lambda}
    \left(
      \frac{\lambda+\bar\lambda}{t_\lambda+2\gamma}
      -
      \frac{\lambda-\bar\lambda}{t_\lambda}
    \right)
    &= -
       \frac{1}{4\pi\gamma}
       \,
       \partial_{\bar\lambda}
       \left(
         (\lambda+\bar\lambda)
         \left(
           1 - \frac{\lambda-\bar\lambda}{\omega-\bar\omega}
         \right)
         -
         (\lambda-\bar\lambda)
         \left(
           \frac{\omega-\bar\omega}{\lambda-\bar\lambda} - 1
         \right)
       \right) \\
    &= -
       \frac{1}{4\pi\gamma}
       \,
       \partial_{\bar\lambda}
       \left(
         (\lambda+\bar\lambda)
         -
         (\omega+\bar\omega)
         -
         (\omega-\bar\omega)
         +
         (\lambda-\bar\lambda)
       \right) \\
    &= -
       \frac{1}{4\pi\gamma}
       \,
       \partial_{\bar\lambda}
       \left(
         2\lambda
         -
         2\omega
       \right) \\
    &= 0
    ;
  \end{align*}
  Thus we get that the density is constant
  $
  \frac{1}{4\pi}
  \left(
    \frac{1}{\alpha}
    +
    \frac{1}{\beta}
  \right)
  $ on the interior of the ellipse~\eqref{eq:sigma(Salpha+iSbeta)}.
%
%
%
  
\end{Example}
The elliptic law appears in the random matrix literature in
 \cite{{Girko:1997:elliptic}}.

\section{Other examples}
There are some other examples that can be done by ad-hoc methods.
\begin{Example}
  Consider two freely independent symmetries $u_2$ and $v_2$ of trace zero,
  for example the generators of the left regular representation of $\IZ_2*\IZ_2$.
  Here we compute the Brown measure of $T=\alpha u_2 + \beta v_2$.
  To get its spectrum, look at its square
  $$
  (\alpha u_2 + \beta v_2)^2 = \alpha^2+\beta^2 + \alpha\beta(u_2v_2+v_2u_2)
  $$
  Since $u_2v_2=(v_2u_2)^*$ is a Haar unitary, we see that $T^2$ 
   is a normal element with spectrum
  $\sigma( T^2 ) = \alpha^2+\beta^2 + \alpha\beta[-2,2]$.
  Since $T$ and $-T$ have the same distribution, it follows that
  $$
  \sigma(\alpha u_2 + \beta v_2)
  = \left\{
      \pm\sqrt{\alpha^2 + \beta^2 + \alpha\beta t} 
      :
      t \in [-2,2]
  \right\}
  $$
  The Brown measure can be deduced by the same symmetry considerations,
  but for the sake of simplicity let us consider the special case $\alpha=1$, 
  $\beta=i$ only.
  Here the spectrum is the union of the complex intervals $[-1-i,1+i]$ and 
  $[-1+i,1-i]$.
  The Brown measure of $(u_2+iv_2)^2=i(u_2v_2+v_2u_2)$ is the arcsine law
  (we are taking the real part of a Haar unitary)
  $$
  d\nu(t) = \frac{dt}{\pi \sqrt{4-t^2}}
  $$
  on the imaginary axis. By symmetry considerations we must have the same measure on each of 
  the four ``legs'' of the spectrum, call it $\mu_0$, which must satisfy
  \begin{align*}
    \int_0^{\sqrt{2}} f(t^2) d\mu_0(t)
    &= \frac{1}{2}\int_0^2 f(t) \frac{dt}{\pi\sqrt{4-t^2}}\\
    &= \int_0^{\sqrt{2}} f(u^2) \frac{u}{\pi\sqrt{4-u^4}}\, du
  \end{align*}
  and it follows that the density of the Brown measure is
  $$
  d\mu
  \left(
    \frac{1\pm i}{\sqrt{2}}\,t
  \right)
  = d\mu_0(|t|) = \frac{|t|}{\pi \sqrt{4-t^4}}\,dt
  $$
\end{Example}

\begin{Example}
  Other examples that are perhaps attackable arise from the following matrix models.
  Consider $U_2+A$, where $U_2\in U(2N)$ is a unitary matrix s.t.\ $U_2=U_2^*$ and
  $\tr U_2=0$, while $A$ is an arbitrary $2N\times 2N$ matrix.
  The spectrum of $U_2+A$ can be bounded as follows.
  Assume $x$ is a unit eigenvector of $U_2+A$ with eigenvalue $\lambda$,
  then it can be decomposed along the spectral projections of $U_2$:
  $x=x_+ + x_-$ 
  so that $U_2 x = x_+ - x_-$.
  By assumption we also have $(U_2+A)(x_+ + x_-) = \lambda (x_+ + x_-)$,
  and thus
  $$
  x_+ = \frac{1}{2}(1+\lambda-A) x
  \qquad
  x_- = \frac{1}{2}(1-\lambda-A) x
  ;
  $$
  now by orthogonality $\langle x_+, x_- \rangle = 0$ we get
  \begin{align*}
    0 &= \langle (1+\lambda -A)x,(1-\lambda+A)x \rangle \\
      &= (1+\lambda)(1-\bar\lambda) \|x\|^2
         + (1+\lambda) \langle x,Ax\rangle 
         - (1-\bar\lambda) \langle Ax,x\rangle 
         - \|Ax\|^2 \\
      &= (1+\lambda-\bar\lambda-|\lambda|^2)\|x\|^2
         + (\lambda+1) \overline{\langle Ax,x\rangle}
         + (\bar\lambda-1) \langle Ax,x\rangle
         - \|Ax\|^2
    .
  \end{align*}
  Separating real and imaginary part results in two equations
  \begin{gather*}
    1-|\lambda|^2 - \|Ax\|^2
    + \lambda\, \overline{\langle Ax,x\rangle}
    + \bar\lambda\, \langle Ax,x\rangle = 0 \\
    \lambda-\bar\lambda + \overline{\langle Ax,x\rangle} - \langle Ax,x\rangle = 0
  \end{gather*}

  Let us no consider two specific cases.
  \begin{description}
   \item[$A$ is unitary] In this case $\|Ax\|=1$ and $\rho=\langle Ax,x \rangle$
    satisfies the following equations
    \begin{gather*}
      -|\lambda|^2 + \lambda\bar\rho +\bar\lambda\rho = 0 \\
      \lambda-\bar\lambda = \rho-\bar\rho
    \end{gather*}
    or in other words
    \begin{gather*}
      |\lambda-\rho|^2 = |\rho|^2 \\
      \Im\lambda = \Im\rho
    \end{gather*}
    thus $\lambda-\rho$ is real and we have
    $$
    \lambda-\rho = \pm|\rho|
    $$
    i.e., 
    $$
    \lambda\in\{ \rho\pm|\rho| : \rho = \langle Ax,x \rangle \in \co\sigma(A) \}
    $$
   \item[$A=iB$ is purely imaginary]
    Here we assume $A+A^*=0$ and the equations are
    \begin{gather*}
      1 - |\lambda|^2-\|Bx\|^2 + i(\bar\lambda-\lambda) \langle Bx,x\rangle = 0 \\
      \lambda-\bar\lambda = 2i\, \langle Bx,x \rangle
    \end{gather*}
    Hence 
    \begin{gather*}
      \Im \lambda = \langle Bx,x \rangle \\
      (\Re \lambda)^2 = 1 - \|Bx\|^2 + \langle Bx,x \rangle^2
    \end{gather*}
  \end{description}

  If one puts $A=U U_3 U^*$, where $U_3$ is an $6N\times6N$ model of the generator of $\IZ_3$,
  and $U$ is a random unitary $6N\times 6N$ matrix,
  then possible eigenvalues are enclosed by the region shown in figure~\ref{fig:U2+U3_spec}.
  \begin{figure}[ht]
    \begin{center} %
     \includegraphics[width=\textwidth, height=.44\textheight, keepaspectratio]{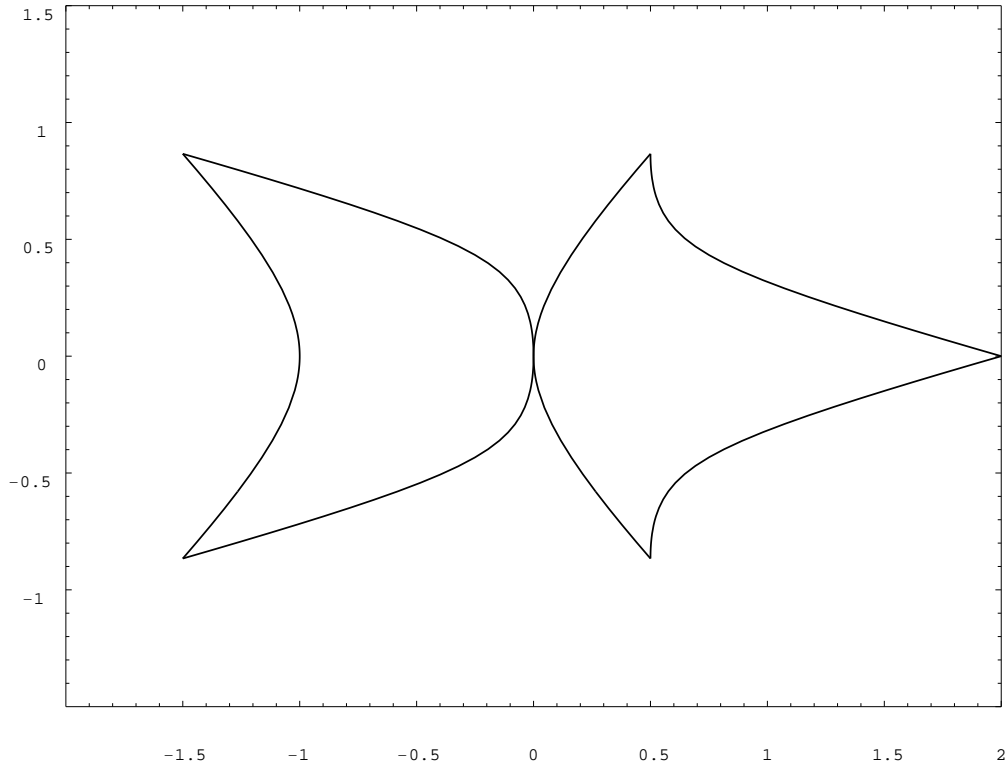}  
      \caption{Possible spectra of random $U_2+U_3$}
      \label{fig:U2+U3_spec}
    \end{center}
  \end{figure}
  And indeed, samples of small numeric random unitary matrices $U_2 + UU_3U^*$ have an
  eigenvalue density as shown in figure~\ref{fig:U2+U3_6},
  while in bigger dimensions the eigenvalues concentrate, cf.\ figure~\ref{fig:U2+U3_150}.
  \begin{figure}[pt]
    \begin{center} %
     \includegraphics[width=\textwidth, height=.44\textheight, keepaspectratio]{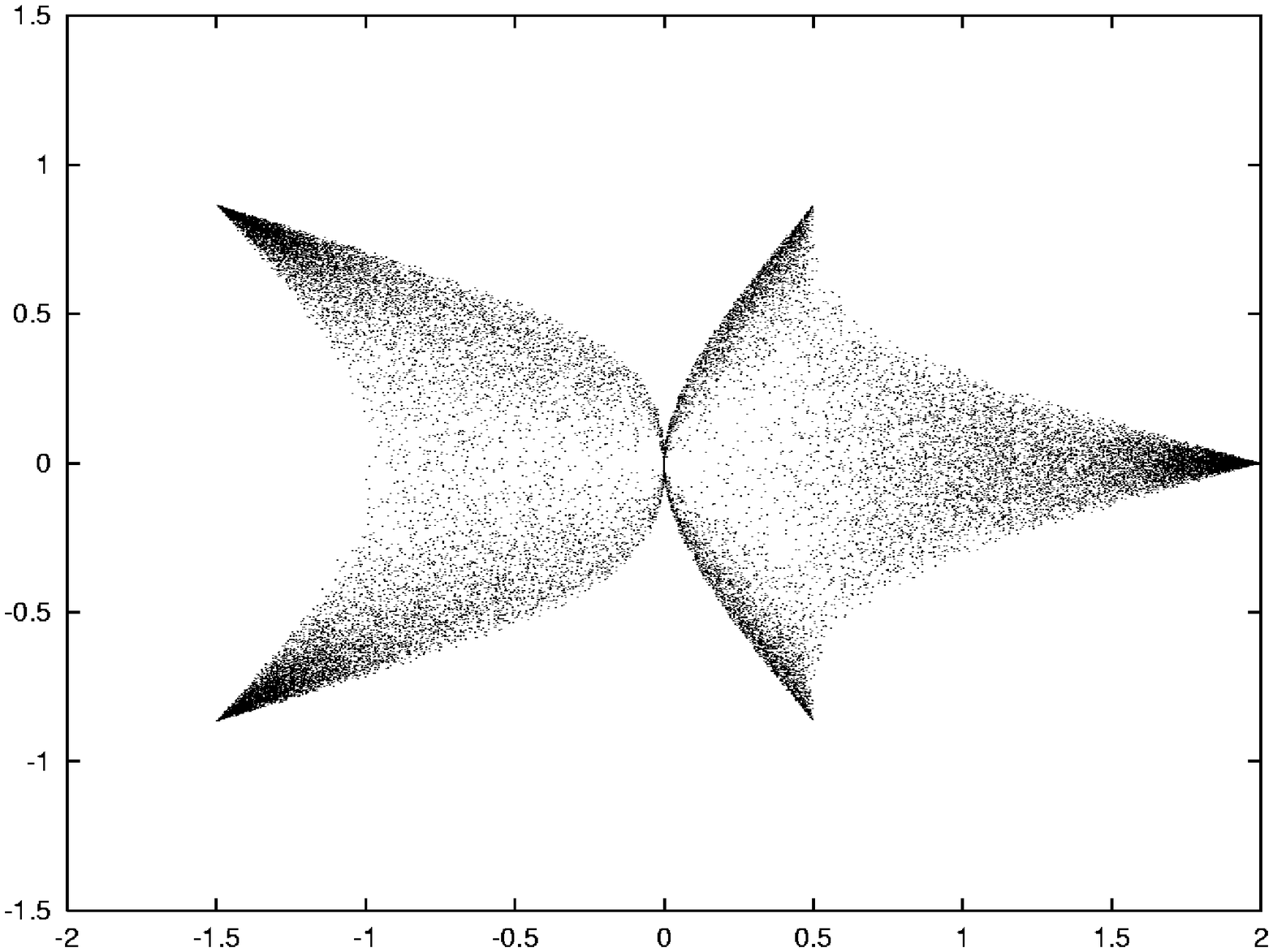}  
      \caption{$5000$ samples of eigenvalues of $6\times6$ random matrices $U_2+U_3$}
      \label{fig:U2+U3_6}
    \end{center}
  \end{figure}
  \begin{figure}[pb]
    \begin{center} %
     \includegraphics[width=\textwidth, height=.44\textheight, keepaspectratio]{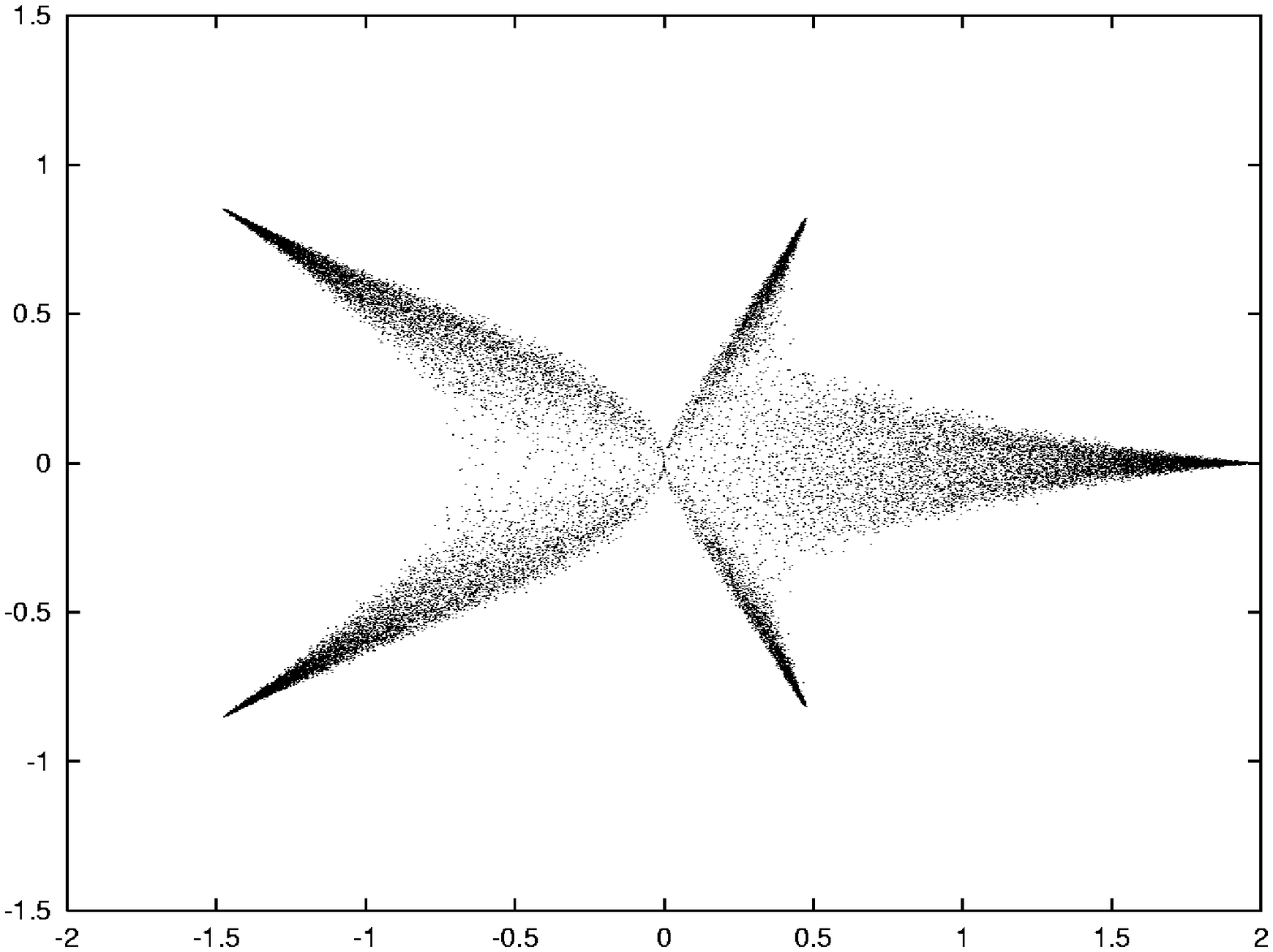}
      \caption{$200$ samples of eigenvalues of $150\times150$ random matrices $U_2+U_3$}
      \label{fig:U2+U3_150}
    \end{center}
  \end{figure}
  We were able to compute the spectrum of the free sum $u_2+u_3$
  recently and will investigate this topic further in future work.

\end{Example}

\bibliography{free,BrownMeasure}
\bibliographystyle{mamsalpha}
\end{document}